\documentclass[10pt]{article}
\usepackage{amssymb,amsfonts,amsmath, amsthm, amsbsy}
\usepackage{color,graphicx, subcaption, placeins, array, mathtools, url, color, rotating} 
\usepackage[text={7in,9.5in},centering,letterpaper]{geometry}
\usepackage[hypertexnames=false,colorlinks=true,urlcolor=blue,linkcolor=blue,citecolor=blue]{hyperref}
\usepackage{tikz}
\usepackage[shortlabels]{enumitem}
\usepackage[margin = 10pt, font=small, labelfont=bf, labelsep = endash]{caption}

\everymath={\displaystyle}
 \numberwithin{equation}{section}

\newtheorem{Theorem}{Theorem}[section]

\newtheorem{Remark}[Theorem]{Remark}
\newtheorem{Assumption}{Assumption}

\def\Re{\mathop{\mathrm{Re}}}

\newcommand{\eps}{\varepsilon}

\title{ A stabilizing effect of advection on planar interfaces in singularly perturbed reaction-diffusion equations}
\author{Paul Carter\thanks{Department of Mathematics, University of California, Irvine, USA}}

\begin{document}
\maketitle

\begin{abstract}
We consider planar traveling fronts between stable steady states in two-component singularly perturbed reaction-diffusion-advection equations, where a small quantity $\delta^2$ represents the ratio of diffusion coefficients. The fronts under consideration are large amplitude and contain a sharp interface, induced by traversing a fast heteroclinic orbit in a suitable slow fast framework. We explore the effect of advection on the spectral stability of the fronts to long wavelength perturbations in two spatial dimensions. We find that for suitably large advection coefficient $\nu$, the fronts are stable to such perturbations, while they can be unstable for smaller values of  $\nu$. In this case, a critical asymptotic scaling $\nu\sim \delta^{-4/3}$ is obtained at which the onset of instability occurs. The results are applied to a family of traveling fronts in a dryland ecosystem model.
\end{abstract}






\section{Introduction}
We consider two-component reaction-diffusion-advection equations of the form
\begin{align}
\begin{split}\label{eq:grda}
U_t &= \Delta U +F(U,V;\boldsymbol{\mu}) \\
V_t &= \frac{1}{\delta^2}\Delta V +G(U,V;\boldsymbol{\mu})+\nu V_x
\end{split}
\end{align}
where $U(x,y,t), V(x,y,t):\mathbb{R}^2\times \mathbb{R}\to \mathbb{R}$, $F$ and $G$ are smooth functions, and $\boldsymbol{\mu}\in\mathbb{R}^m$ denotes a collection of system parameters. We assume that~\eqref{eq:grda} is singularly perturbed, with $0<\delta\ll 1$. The advection coefficient $0\leq \nu<\infty$ is arbitrary. We consider planar interfaces between spatially homogeneous stable steady state solutions $(U,V)(x,y,t) = (U^\pm, V^\pm)$ of~\eqref{eq:grda}. The planar interfaces manifest as traveling wave solutions $(U,V)(x,y,t) = (u_\mathrm{h},v_\mathrm{h})(\xi), \xi=x-ct$, which propagate with constant speed $c$ in the $x$-direction, and are constant in the $y$-direction, and asymptotically approach the steady states $\lim_{\xi\to\pm\infty}(u_\mathrm{h},v_\mathrm{h})(\xi)=(U^\pm, V^\pm)$. 

Reaction-diffusion-advection systems arise in models of diverse phenomena such as pattern formation in mussel beds~\cite{bennett2019large} and plankton~\cite{malchow1996nonlinear}, fog and wind induced vegetation alignment~\cite{borthagaray2010vegetation}, disease spread~\cite{ge2015sis}, and population dynamics~\cite{chen2008evolution}. Here, we are primarily motivated by the phenomenon of desertification fronts in water-limited ecosystems~\cite{zelnik2017desertification}, in which the bare-soil state slowly invades a vegetated state, resulting in (typically irreversible) desertification~\cite{hellden1988desertification, un2015transforming}, and similarly the reverse mechanism of vegetation fronts, in which vegetation invades a bare soil state. Instabilities in the resulting planar interface between vegetation and bare soil have been linked to spatial pattern formation~\cite{banerjee2023rethinking,CDLOR,fernandez2019front}; see also~\S\ref{sec:discussion}. Similar interfaces also appear in savanna/forest ecosystems, cloud formation, salt marshes, and other applications~\cite{bastiaansen2022fragmented, lipcius2021facilitation}.

Our aim is to examine the effect of the advection term $\nu V_x$ on the stability of such an interface in two spatial dimensions. In systems of the form~\eqref{eq:grda}, the relation between diffusive and advective dynamics is known to impact the stability of planar stripes and periodic patterns~\cite{siero2015striped,siteur2014beyond}, and in particular the presence of advection can have a stabilizing effect in the direction along the stripe~\cite{kolokolnikov2018stabilizing}. We aim to explore the stabilizing effect of advection, focussing on long wavelength instabilities in the case of a single planar interface between stable steady states. We note that such interfaces in the class of equations~\eqref{eq:grda} was previously studied in the absence of advection, i.e. $\nu=0$, on infinite~\cite{CDLOR} and cylindrical domains~\cite{taniguchi2003instability, taniguchi1994instability}.  

In the context of the motivating example of vegetation pattern formation in dryland ecosystem models, the quantities $U$ and $V$ represent interacting species and/or resources, for example $U$ may represent vegetation biomass, and $V$ represents water availability. In such models, it is natural to have widely separated diffusion coefficients due to the differing length and/or time scales on which water is transported and on which different vegetation species evolve~\cite{rietkerk2008regular}. In this setting, the advection term represents a slope in the topography, leading to downhill flow of water, and thus an anisotropy in the system. Observations suggest that the absence of advection (that is, flat terrain) lead to spotted and/or labrynthine patterns, whereas on sloped terrain vegetation may align in bands, consisting of interfaces alternating between vegetated and desert states~\cite{barbier2014case,deblauwe2012determinants, deblauwe2011environmental, gandhi2018influence, ludwig2005vegetation,valentin1999soil}. These interfaces align perpendicular to the slope, suggesting that the downhill flow of water prescribes a preferred orientation of the interface~\cite{meron2012pattern}.  In~\cite{CDLOR} it was shown that such planar interfaces are unstable in many ecosystem models in the absence of advection, and it is the goal of this work to examine the effect of advection on the (in)stability of such interfaces. In particular, we demonstrate that sufficiently large advection has a stabilizing effect on interfaces in two spatial dimensions, with respect to long wavelength perturbations in the $y$-direction, transverse to the direction of propagation.

In the spirit of~\cite{CDLOR}, our results are framed in the context of a geometric singular perturbation analysis of the traveling wave equation associated with~\eqref{eq:grda}, under suitable assumptions about the underlying geometry of the system. The novel contribution of the current study is the inclusion of the advection term; the coefficient $\nu$ can be small or large relative to the small parameter $\delta$, which naturally leads to a three timescale system, which must be separately analyzed in several scaling regimes. Nevertheless, we obtain a simple explicit criterion for (in)stability of planar interfaces, depending on the relative scaling of $\nu, \delta$, through a formal asymptotic approach, and we identify a potential onset of (in)stability at a critical scaling $\nu\sim \delta^{-4/3}$. The results can be easily applied to fronts in reaction-diffusion-advection models, and we demonstrate the applicability of the results to a dryland ecosystem model in~\S\ref{sec:dryland_example}.

We note that the results also apply to systems of the form
\begin{align}
\begin{split}
U_t &= \Delta U +F(U,V;\boldsymbol{\mu}) +\nu_1 U_x\\
V_t &= \frac{1}{\delta^2}\Delta V +G(U,V;\boldsymbol{\mu})+\nu_2 V_x
\end{split}
\end{align}
for arbitrary advection coefficients $\nu_i\in \mathbb{R}$. By shifting to a traveling coordinate frame, and reversing the spatial variable $x$ if necessary, this system can be transformed to~\eqref{eq:grda}, defining $0\leq \nu =|\nu_1-\nu_2|<\infty$ as the differential flow~\cite{rovinsky1992chemical,siero2015striped}. Additionally, we note that the reaction terms $F$ and $G$ in~\eqref{eq:grda} do not depend explicitly on $\delta$ or $\nu$. While one could consider such a dependence in a given model with explicit reaction terms, we will see that the behavior of this system depends critically on certain relative scalings between the parameters $\delta$ and $\nu$. To avoid additionally tracking these scalings within the reaction terms themselves, for simplicity we assume they are independent of $(\delta,\nu)$.

\section{Setup}
\subsection{Traveling wave formulation}\label{sec:travelingwave_setup}
To capture traveling front solutions which propagate in the direction determined by the advection term, we move into a traveling coordinate frame $\xi=x-ct$ and obtain the system
\begin{align}
\begin{split}\label{eq:grda_tf}
U_t &=  U_{\xi\xi}+U_{yy}+cU_\xi +F(U,V) \\
V_t &= \frac{1}{\delta^2}(V_{\xi\xi}+V_{yy}) +\left(\nu+c\right) V_\xi+G(U,V)
\end{split}
\end{align}
where we drop the explicit dependence on the system parameters $\boldsymbol{\mu}$. We search for stationary solutions $(u,v)(\xi,y,t) = (u,v)(\xi)$, which are constant in the $y$-direction, and thus propagate with constant speed $c$ in the $x$-direction. This results in the traveling wave ODE
\begin{align}
\begin{split}\label{eq:TW_ode}
0&= u_{\xi\xi}+cu_\xi+F(u,v)\\
0&=v_{\xi\xi}+\delta^2(\nu+c)v_\xi +\delta^2 G(u,v).
\end{split}
\end{align}
The system~\eqref{eq:TW_ode} can then be written as the first order system
\begin{align}
\begin{split}\label{eq:fast}
u_\xi&= p\\
 p_\xi&= -cp-F(u,v)\\
v_\xi&= \delta q\\
q_\xi&= -\delta^2(\nu+c)q-\delta G(u,v).
\end{split}
\end{align}
where the homogeneous rest states $(U^\pm,V^\pm)$ of~\eqref{eq:grda} are given by the fixed points $P^\pm = (U^\pm,0,V^\pm,0)$ of~\eqref{eq:fast}.

To analyze traveling front solutions in~\eqref{eq:TW_ode}, we use geometric singular perturbation theory~\cite{fenichel1979geometric}. Throughout, we assume that $\delta\ll1$ is a small parameter, but the parameter $\nu$ can be small or large.  Thus in the regime $\delta \ll1$,  this system can have up to three timescales, determined by the relation between the two parameters $\nu,\delta$, and we must therefore separate the analysis of~\eqref{eq:TW_ode} into cases, depending on the relative size of the parameters $\nu,\delta$. 

\begin{figure}[t]
\centering
\includegraphics[width=0.4\linewidth]{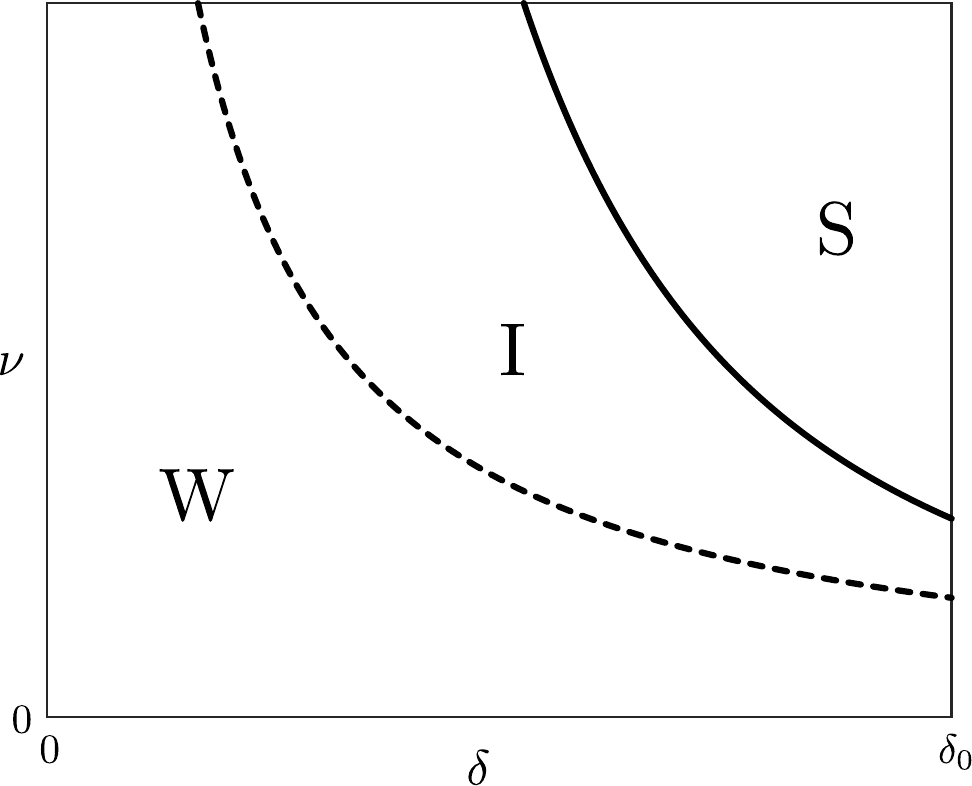}
\caption{Shown is a schematic of the different scaling regimes considered. The dashed and solid curves represent the boundary curves $\nu = \tfrac{1}{\delta}$ and $ \nu = \tfrac{r_0}{ \delta^2}$, respectively, as described in~\S\ref{sec:summary}. The weak advection (W), intermediate (I), and strong advection (S) regimes are labeled accordingly.  }
\label{fig:scaling_regimes}
\end{figure}

The following singular perturbation analysis distinguishes between $3$ cases, depending on which parameter is used as the primary singular perturbation parameter, and we describe the slow/fast structure of traveling fronts in each of these regions for sufficiently small $0<\delta \ll 1$. In the weak advection regime $0\leq \nu \leq \mathcal{O}(\delta^{-1})$, $\delta$ serves as the timescale separation parameter, while in the strong advection regime $\nu \geq \mathcal{O}(\delta^{-2})$, the quantity $\eps:=\nu^{-1}\ll 1$ is taken as the timescale separation parameter. In the intermediate regime $\mathcal{O}(\delta^{-1}) \leq  \nu \leq \mathcal{O}(\delta^{-2})$, the advection-diffusion coefficient ``ratio" $r := \delta^2 \nu$ is taken as the primary singular perturbation parameter. By combining the results in these regions, we are able to describe the slow/fast structure of traveling front solutions for each $(\nu,\delta)$ satisfying $\nu \geq0, 0<\delta\leq \delta_0$ for some $\delta_0>0$; see Figure~\ref{fig:scaling_regimes} and~\S\ref{sec:summary}. Due to the appearance of the $\delta^2\nu$ coefficient in~\eqref{eq:TW_ode}, the quantity $r$ will in fact play an important role throughout the three regimes.


In each regime, the set $\mathcal{S}=\{F(u,v)=0\}$ organizes the dynamics, as this set helps define the critical manifold(s) which appear in the slow-fast formulation of the traveling wave problem.  Generically, on a given compact set, away from points where $F_u=0$, $\mathcal{S}$ is formed by the union of a finite number of branches $\mathcal{S}_j$, $j=1,\ldots,N$, which can be written as graphs $u=f_j(v)$ satisfying $F(f_j(v),v))=0$ for $v\in I_j$, where $I_j$ is an interval. We denote by $f^\pm$ the functions which define the graphs corresponding to the two branches $\mathcal{S}^\pm$ of $\mathcal{S}$ satisfying $f^\pm(U^\pm)=V^\pm$, and we let $F_u^\pm$ denote $\tfrac{\partial F}{\partial u}(U^\pm, V^\pm)$, etc. Independent of the specific parameter regime, we make the following basic assumptions regarding the steady states $(U^\pm,V^\pm)$.

\begin{Assumption}\label{assump:steadystates} (Steady states)
\begin{enumerate}[(i)]
    \item \label{assump:bistable} There exist two homogeneous steady states $(U^\pm,V^\pm)$ which are stable as solutions of~\eqref{eq:grda} for $0<\delta\ll1$ and all $\nu\geq 0$.  In particular, we assume (see Appendix~\ref{app:steadystates})
    \begin{align}\label{eq:steadystate_conditions}
F_u^\pm<0,\qquad  G_v^\pm<0, \qquad F_u^\pm G_v^\pm- F_v^\pm G_u^\pm>0.
    \end{align}
     \item \label{assump:branches}  The states $(U^\pm,V^\pm)$ lie on different branches of $\mathcal{S}$, that is, $f^-(v)\not\equiv f^+(v)$.
    \end{enumerate}
\end{Assumption}

\begin{figure}
\hspace{.05\textwidth}
\begin{subfigure}{.35 \textwidth}
\centering
\includegraphics[width=1\linewidth]{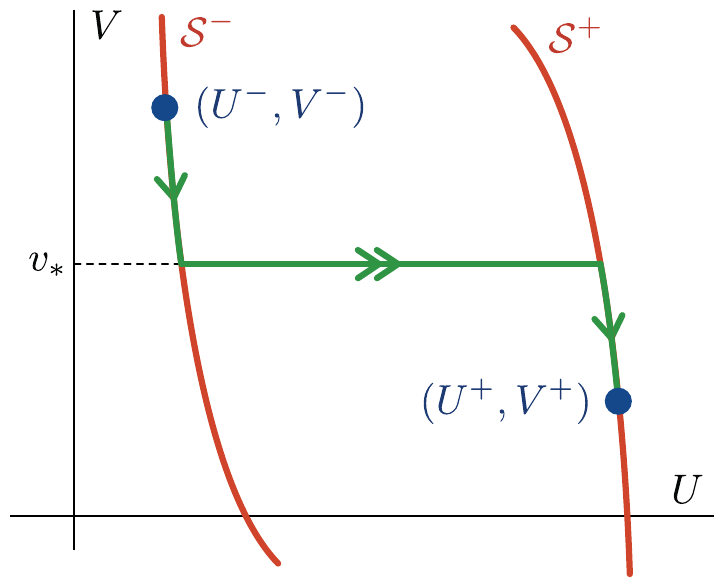}
\end{subfigure}
\hspace{.1\textwidth}
\begin{subfigure}{.4 \textwidth}
\centering
\includegraphics[width=1\linewidth]{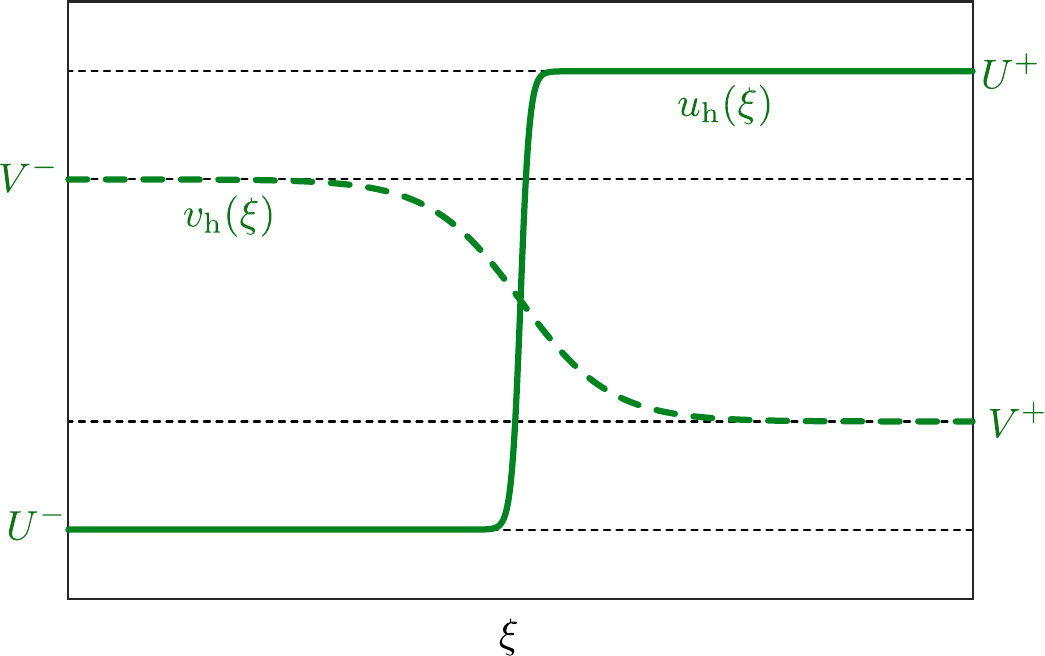}
\end{subfigure}
\caption{(Left) The steady states $(U^\pm,V^\pm)$ and the branches $\mathcal{S}^\pm$ of the nullcline $F(U,V)=0$.  (Right) Schematic of a traveling front solution $(u_\mathrm{h},v_\mathrm{h})(\xi)$ with a single sharp interface.   }
\label{fig:schematic_front}
\end{figure}

We will analyze traveling front solutions $(u,v)(\xi)$ satisfying $\lim_{\xi \to \pm \infty}(u,v)(\xi) = (U^\pm, V^\pm)$. The first condition~\ref{assump:bistable} ensures that such a front is \emph{bistable}, so that it forms an interface between asymptotically stable rest states; we will see that the resulting conditions~\eqref{eq:steadystate_conditions} ensure hyperbolicity of relevant critical manifolds and the rest states in their associated reduced flows. We also remark that the condition on $G^\pm_v$ ensures the steady states remain stable for large advection $\nu\gg1$ (see Appendix~\ref{app:steadystates}); this prevents instabilities which can arise in the background states $(U^\pm, V^\pm)$ due to large differential flow~\cite{borgogno2009mathematical, rovinsky1992chemical, siero2015striped}, so that we focus only on instabilities which arise due to the interface itself.

The second condition~\ref{assump:branches} ensures that the front interface is sharp; that is, in the appropriate slow/fast formulation (depending on the specific asymptotic regime of the parameters $\delta, \nu$), the front (in an appropriate singular limit) must traverse a singular fast heteroclinic orbit $u_*(\xi)$ of an associated layer problem in the subspace $v=v_*$ with leading order speed $c=c_*$, as opposed to being entirely contained within the reduced flow on a single connected branch of a slow manifold; for the latter situation, see e.g.~\cite{doelman2022slow}. The solution $u_*(\xi)$ solves the simpler scalar traveling wave equation
\begin{align}\label{eq:reduced_fast}
0&= u_{\xi\xi}+c_*u_\xi+F(u,v_*),
\end{align}
and forms a fast connection between the branches of $\mathcal{S}$, that is, $u_*^\pm:=\lim_{\xi\to\pm\infty}u_*(\xi)$ satisfy $u_*^\pm = f^\pm(v_*)$. The existence of such a heteroclinic orbit in~\eqref{eq:reduced_fast} can for instance be obtained in a given system using phase plane techniques. The fronts under consideration here traverse only one such fast heteroclinic orbit, i.e. they do not jump back and forth between different branches of $\mathcal{S}$.

Beyond the conditions~\ref{assump:bistable}-\ref{assump:branches} on the steady states $(U^\pm,V^\pm)$, additional structure is required concerning the reduced flows of certain critical manifolds which appear in the existence problem to allow for bistable traveling front solutions between $(U^\pm, V^\pm)$. However, since these conditions are related to the specific slow/fast formulation in each parameter regime, we delay their discussion until introducing the slow/fast structure of the fronts; see Assumptions~\ref{assump:diffusion_existence} and~\ref{assump:advection_existence} in~\S\ref{sec:diffusion_dominant} and~\S\ref{sec:advection_dominant}, respectively. In the following section, we focus on the stability criterion which arises given a traveling wave solution which connects the steady states $(U^\pm, V^\pm)$.

\subsection{Long-wave (in)stability of planar interfaces}

Given a traveling wave solution $\phi_\mathrm{h}(\xi;r, \delta) = (u_\mathrm{h},v_\mathrm{h})(\xi;\nu, \delta)$ of~\eqref{eq:grda_tf} with speed $c=c(\nu, \delta)$ satisfying $\lim_{\xi\to\pm\infty}\phi_\mathrm{h}(\xi;\nu, \delta)=(U^\pm, V^\pm)$, we have the corresponding linear stability problem
\begin{align}
\begin{split}\label{eq:mKstability}
\lambda u &=u_{\xi\xi}+cu_\xi -\ell^2u +F_u(u_\mathrm{h}(\xi),v_\mathrm{h}(\xi))u+F_v(u_\mathrm{h}(\xi),v_\mathrm{h}(\xi))v\\
\lambda v &= \frac{1}{\delta^2}v_{\xi\xi}+\left(\nu+c\right)v_\xi-\frac{\ell^2}{\delta^2}v+G_u(u_\mathrm{h}(\xi),v_\mathrm{h}(\xi))u+G_v(u_\mathrm{h}(\xi),v_\mathrm{h}(\xi))v,
\end{split}
\end{align}
parameterized by the transverse wavenumber $\ell\in\mathbb{R}$, which can equivalently be written as
\begin{align}\label{eq:eigenvalue_problem}
\mathcal{L}\begin{pmatrix} u\\ v \end{pmatrix} =\lambda \begin{pmatrix} u\\ v \end{pmatrix} + \ell^2 \begin{pmatrix} u\\ \frac{1}{\delta^2}v \end{pmatrix} .
\end{align}
where 
\begin{align}
\mathcal{L}=\begin{pmatrix}\partial_{\xi\xi}+c\partial_\xi +F_u(u_\mathrm{h}(\xi),v_\mathrm{h}(\xi)) & F_v(u_\mathrm{h}(\xi),v_\mathrm{h}(\xi))\\ G_u(u_\mathrm{h}(\xi),v_\mathrm{h}(\xi))& \frac{1}{\delta^2}\partial_{\xi\xi}+\left(\nu+c\right)\partial_\xi+ G_v(u_\mathrm{h}(\xi),v_\mathrm{h}(\xi))\end{pmatrix}.
\end{align}

When $\ell=0$, this eigenvalue problem is solved by taking $\lambda=0$ and $(u,v) = (u_\mathrm{h}',v_\mathrm{h}')$ (due to translation invariance). We make the following assumption regarding the stability of the front as a traveling wave in one space dimension, i.e. in the direction of propagation.

\begin{Assumption} (1D stability of the front) 
    The operator $\mathcal{L}$ satisfies $\mathrm{spec}\{\mathcal{L}\} \subset \{ \lambda \in \mathbb{C}: \Re \lambda<0\} \cup \{0\}$. Furthermore, the eigenvalue $\lambda=0$ is isolated and algebraically simple, and $\mathcal{L}$ has one-dimensional generalized kernel spanned by the eigenfunction $(u_\mathrm{h}',v_\mathrm{h}')$.
\end{Assumption}

To study the stability of the front to long wavelength perturbations in two spatial dimensions, we examine how this eigenvalue problem perturbs for values of $|\ell|\ll 1$. Following the  (formal) analysis of~\cite{CDLOR} for the stability problem~\eqref{eq:mKstability}, we expand the critical translation eigenvalue $\lambda_\mathrm{c}(\ell)$ satisfying $\lambda_\mathrm{c}(\ell)=0$ and the corresponding eigenfunction $(u_\mathrm{c},v_\mathrm{c})(\xi;\ell)$ as
\begin{align}
\lambda_\mathrm{c}(\ell)&= \lambda_{\mathrm{c},2} \ell^2+\mathcal{O}(\ell^4), \qquad \begin{pmatrix} u_\mathrm{c}(\xi;\ell)\\ v_\mathrm{c}(\xi;\ell)\end{pmatrix}=\begin{pmatrix} u_\mathrm{h}'(\xi)\\ v_\mathrm{h}'(\xi)\end{pmatrix}+\ell^2\begin{pmatrix} u_\mathrm{c,2}(\xi)\\ v_\mathrm{c,2}(\xi)\end{pmatrix}+\mathcal{O}(\ell^4)
\end{align}
Substituting into~\eqref{eq:eigenvalue_problem}, at $\mathcal{O}(\ell^2)$, we have the equation
\begin{align}
\mathcal{L}\begin{pmatrix} u_\mathrm{c,2}\\ v_\mathrm{c,2}\end{pmatrix} =\lambda_{\mathrm{c},2} \begin{pmatrix} u_\mathrm{h}'\\ v_\mathrm{h}'\end{pmatrix} +  \begin{pmatrix} u_\mathrm{h}'\\ \frac{1}{\delta^2}v_\mathrm{h}' \end{pmatrix},
\end{align}
which leads to the Fredholm solvability condition
 \begin{align}\label{eq:solvability}
\left\langle \lambda_{\mathrm{c},2} \begin{pmatrix} u_\mathrm{h}'\\ v_\mathrm{h}'\end{pmatrix} +  \begin{pmatrix} u_\mathrm{h}'\\ \frac{1}{\delta^2}v_\mathrm{h}' \end{pmatrix}, \begin{pmatrix} u^A\\ v^A\end{pmatrix} \right\rangle_{L^2}=0
\end{align}
where $(u^A,v^A)(\xi)$ is the unique (up to scalar multiple) integrable eigenfunction of the adjoint equation
\begin{align}
\mathcal{L}^A\begin{pmatrix} u\\ v \end{pmatrix} =0,
\end{align}
where the adjoint operator $\mathcal{L}^A$ is given by
\begin{align}
\mathcal{L}^A=\begin{pmatrix}\partial_{\xi\xi}-c\partial_\xi +F_u(u_\mathrm{h}(\xi),v_\mathrm{h}(\xi)) &G_u(u_\mathrm{h}(\xi),v_\mathrm{h}(\xi)) \\ F_v(u_\mathrm{h}(\xi),v_\mathrm{h}(\xi))& \frac{1}{\delta^2}\partial_{\xi\xi}-\left(\nu+c\right)\partial_\xi+ G_v(u_\mathrm{h}(\xi),v_\mathrm{h}(\xi))\end{pmatrix}.
\end{align}
Equivalently, solving~\eqref{eq:solvability} for $\lambda_{\mathrm{c},2}$, we obtain
 \begin{align}\label{eq:lambda2c_expression}
\lambda_{\mathrm{c},2} =-\frac{\int_{-\infty}^\infty \left(u_\mathrm{h}'(\xi)u^A(\xi)+\frac{1}{\delta^2}v_\mathrm{h}'(\xi)v^A(\xi)\right)\mathrm{d}\xi}{\int_{-\infty}^\infty \left(u_\mathrm{h}'(\xi)u^A(\xi)+v_\mathrm{h}'(\xi)v^A(\xi)\right)\mathrm{d}\xi}.
\end{align}
The sign of $\lambda_{\mathrm{c},2}$ determines the $2$D stability of the front $(u_\mathrm{h},v_\mathrm{h})$ to perturbations with small transverse wavenumber $|\ell|\ll1$. As with the slow fast structure of the fronts themselves, the structure of the stability problem~\eqref{eq:mKstability} and the computation of the adjoint solution $(u^A,v^A)(\xi)$ change depending on the relative size(s) of the parameters $\nu,\delta$. Hence we must split the computation of $\lambda_{\mathrm{c},2} $ into cases corresponding to different scaling regimes as described in~\S\ref{sec:travelingwave_setup}.

\subsection{Summary of results}\label{sec:summary}
By considering the slow fast construction of traveling fronts in the weak advection, intermediate, and strong advection regimes, we determine leading order asymptotics for the critical coefficient $\lambda_{\mathrm{c},2}$~\eqref{eq:lambda2c_expression} which determines long wavelength instabilities along the front interface. We will see that the sign of this coefficient depends only on information encoded in the fast layer orbit $u_*(\xi)$ of~\eqref{eq:reduced_fast} in the subspace $v=v_*$ in the singular slow/fast framework. We impose one additional nondegeneracy assumption
\begin{Assumption}\label{assump:nondegeneracy}
    (Nondegeneracy condition) The quantity $G(u^+_*,v_*)-G(u^-_*,v_*)\neq 0$.
\end{Assumption}

Under Assumptions~\ref{assump:steadystates}--\ref{assump:nondegeneracy}, for a traveling front $\phi_\mathrm{h}$ which traverses a single fast jump $u_*(\xi)$ of the reduced equation~\eqref{eq:reduced_fast}, we obtain an asymptotic long-wavelength stability criterion which holds throughout the weak advection, strong advection, and intermediate regimes. To summarize, letting $r_0$ denote a sufficiently small fixed positive constant, for sufficiently small $\delta_0>0$ we find the following asymptotic stability criteria by determining the sign of $\lambda_{\mathrm{c},2}$ for $0<\delta \ll \delta_0 \ll 1$:
\begin{itemize}
    \item \textbf{Weak advection regime:} $0\leq \nu \leq \frac{1}{\delta}$. In this regime, to leading order
    \begin{align}\label{eq:diffusion_criterion}
        \text{sign}(\lambda_{\mathrm{c},2})=-\text{sign}(F_*)\times \text{sign}(G_*)
    \end{align}
    where 
      \begin{align}\label{eq:fstar_gstar}
        F_*: = \int_{-\infty}^\infty F_v(u_*(\xi),v_*)e^{c_*\xi}u_*'(\xi)\mathrm{d}\xi, \qquad G_*:= G(u^+_*,v_*)-G(u^-_*,v_*)
    \end{align}
and $u^\pm_*:=\lim_{\xi\to\pm\infty}u_*(\xi)$.
    \item \textbf{Intermediate regime:} $\frac{1}{\delta} \leq  \nu \leq \frac{r_0}{ \delta^2}$. In this regime, to leading order
     \begin{align}\label{eq:intermediate_criterion}
        \text{sign}(\lambda_{\mathrm{c},2})=\text{sign}\left(-1+ \frac{M}{\delta^4\nu^3} \right)
    \end{align}
    where $M$ = $\mathcal{O}(1)$ with respect to $\delta, \nu$, and $\text{sign}(M)=-\text{sign}(F_*)\times \text{sign}(G_*)$. In particular, if $M>0$, then to leading order, $\lambda_{\mathrm{c},2}$ changes sign when $\nu\sim M^{1/3}\delta^{-4/3}$.
    \item \textbf{Strong advection regime:} $\nu \geq \frac{r_0}{\delta^2}$. Throughout this regime, to leading order we find that
\begin{align}\label{eq:advection_criterion}
\text{sign}(\lambda_{\mathrm{c},2})=-1. \end{align}
\end{itemize}

The asymptotic estimates are uniform in sufficiently small $\delta$, so that the three regimes collectively describe the parameter region $\{(\delta,\nu): 0<\delta\leq \delta_0, 0\leq \nu <\infty\}$ for some suitably small choice of $\delta_0>0$. In the weak advection regime, the stability of traveling fronts to long wavelength perturbations is encoded purely in the nonlinearities $F$ and $G$, evaluated along the fast heteroclinic orbit $u_*(\xi)$. This criterion is analogous to that obtained in~\cite{CDLOR} in the absence of advection, and in the limit $\nu\to0$ the corresponding expression for $\lambda_{\mathrm{c,2}}$ agrees with that found in~\cite[\S2]{CDLOR}. As $\nu$ increases relative to $\delta$, in the intermediate regime, depending on the nonlinearities $F$ and $G$ and the relative size of $\nu$ with respect to $\delta$, traveling fronts can be stable or unstable to long wavelength perturbations, with a potential sign change of $\lambda_{\mathrm{c},2}$ occurring at the critical scaling $\nu \sim \mathcal{O}(\delta^{-4/3})$. Finally, in the strong advection regime, all bistable traveling fronts considered here are stable to long wavelength perturbations. In this sense, the presence of advection has a stabilizing effect on the front as a planar interface. 

To obtain the stability criteria above, we employ a mixture of geometric singular perturbation theory and formal asymptotic arguments to construct the adjoint solution $(u^A, v^A)(\xi)$ and estimate the expression~\eqref{eq:lambda2c_expression} in each of the scaling regimes. However, we emphasize that the results above could in principle be obtained rigorously using geometric singular perturbation methods, in combination with exponential dichotomies/trichotomies and Lin's method, or Evans function approaches; see e.g.~\cite{BCD, sewaltspatially} for thorough analyses of  stability of planar traveling fronts and stripe solutions in specific reaction-diffusion-advection equations. However, for our purposes, we believe that such a technical analysis would detract from the simple message herein, that the presence of advection has a stabilizing effect on planar interfaces, and a straightforward stability criterion which can easily be applied in many example systems.

The remainder of the paper is outlined as follows. In~\S\ref{sec:diffusion_dominant} we describe the construction of traveling fronts and the leading order computation of $\lambda_{\mathrm{c},2}$ in the weak advection and intermediate regimes, while the strong advection regime is considered in~\S\ref{sec:advection_dominant}. In~\S\ref{sec:dryland_example}, we apply these results to an explicit dryland ecosystem model, and~\S\ref{sec:discussion} contains some numerical simulations and a brief discussion of the results.

\section{Weak advection and intermediate regimes}\label{sec:diffusion_dominant}
Due to the similarity in the slow/fast geometry associated with the weak advection and intermediate regimes, in this section we consider both regimes, and outline the differences in each case. Taken together, we consider the regime $0\leq \nu \leq \tfrac{r_0}{\delta^2}$, or equivalently $0\leq r \leq r_0$, where $r_0$ is a (yet to be fixed) small parameter. Therefore, we are interested in the behavior of the traveling wave equation~\eqref{eq:TW_ode} when both $\delta$ and $r=\delta^2\nu$ can be taken as small parameters.

This leads to a system with (up to) three timescales, and there is a distinction between the singular limits obtained by taking $\delta\to0$ with $r$ fixed, versus $r\to0$ with $\delta$ fixed. Hence the case $0\leq r\leq r_0$ needs to be split into two subcases: (i) $r= \bar{r}\delta, \bar{r}\leq \bar{r}_0$, corresponding to the weak advection regime and (ii) $\delta= \bar{\delta}r, \bar{\delta}\leq \bar{\delta}_0$, corresponding to the intermediate regime. We can choose the quantities $\bar{r}_0,\bar{\delta}_0>1$ so that these regimes overlap, and we can understand the slow/fast structure of traveling fronts in the entire region $0\leq r\leq r_0, 0<\delta\leq \delta_0$ for some small $0<r_0,\delta_0\ll1$.

We begin by describing the slow/fast structure of traveling fronts in each case in~\S\ref{sec:existence_diffusion_dominant}, followed by a computation of the coefficient $\lambda_{\mathrm{c},2}$ in~\S\ref{sec:stability_diffusion_dominant}.

\subsection{Structure of traveling fronts}\label{sec:existence_diffusion_dominant}

The structure of the orbits in each case is similar, but with a different parameter used as the timescale separation parameter ($\delta$ vs. $r$) in each case. We consider the first case in detail, and then outline differences relevant for the analysis of the second case. 

\subsubsection{Case $(i)$: $r= \bar{r}\delta, 0\leq \bar{r}\leq \bar{r}_0$}\label{sec:diffusion_existence_i}
We set $r=\bar{r}\delta$ in~\eqref{eq:fast}, obtaining
\begin{align}
\begin{split}\label{eq:fastr0}
u_\xi&= p\\
 p_\xi&= -cp-F(u,v)\\
v_\xi&= \delta q\\
q_\xi&= -\delta\left((\bar{r}+\delta c)q+ G(u,v)\right).
\end{split}
\end{align}
which we aim to analyze for all $0\leq\bar{r}\leq \bar{r}_0$ for fixed $\bar{r}_0>0$, where $\bar{r}_0$ is $\mathcal{O}(1)$ with respect to $\delta$.

\begin{figure}[t]
\centering
\includegraphics[width=0.8\linewidth]{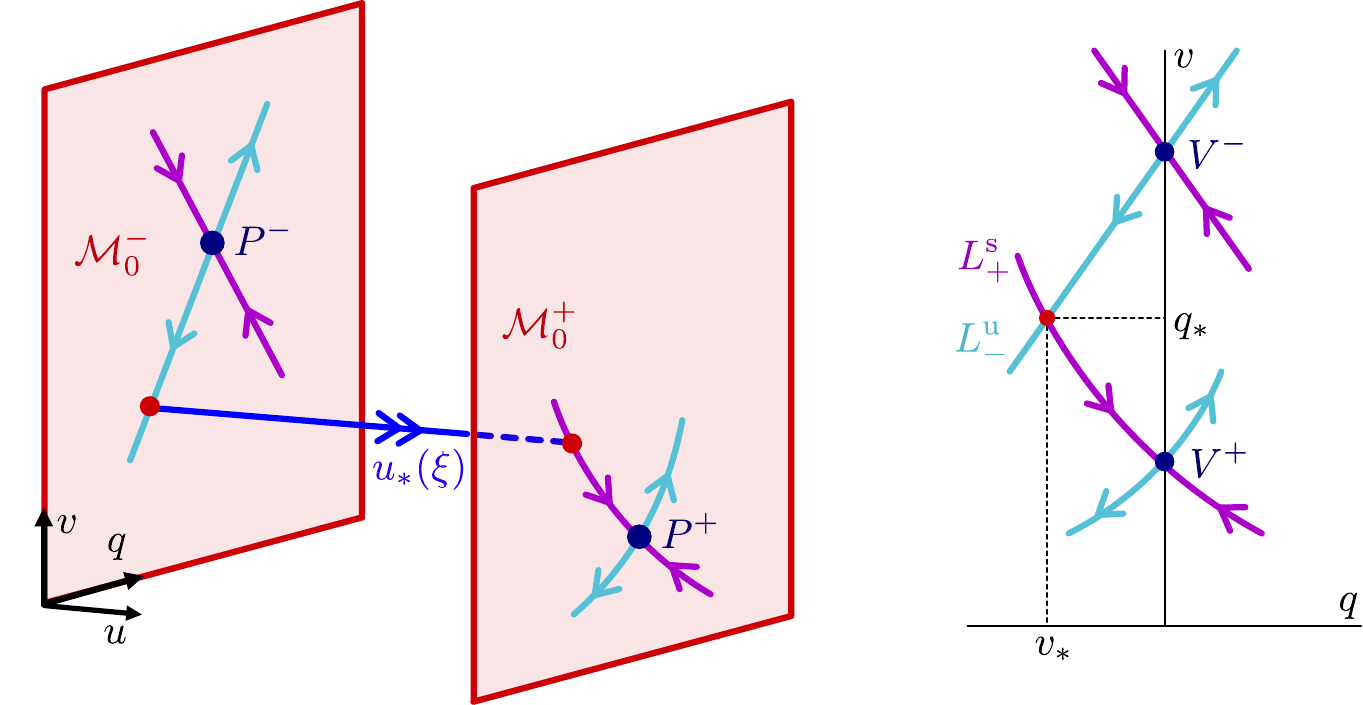}
\caption{Shown is a schematic of the slow/fast construction of the singular traveling front in the weak advection regime. (Left) The fast jump $u_*(\xi)$ of the layer problem~\eqref{eq:layer} between the manifolds $\mathcal{M}^\pm_0$ in $(u,v,q)$-space. (Right) Intersection of (projection of) manifolds $L^\mathrm{u}_-$ and $L^\mathrm{s}_+$ of the reduced flows~\eqref{eq:reduced} as in Assumption~\ref{assump:diffusion_existence}.  }
\label{fig:diffusion_singular}
\end{figure}

This results in a $2$-fast-$2$-slow system with timescale separation parameter $\delta$. Setting $\delta=0$, we obtain the layer problem
\begin{align}
\begin{split}\label{eq:layer}
u_\xi&= p\\
 p_\xi&= -cp-F(u,v)
\end{split}
\end{align}
which, by Assumption~\ref{assump:steadystates}\ref{assump:branches}, admits two hyperbolic equilibria given by $u=f^\pm(v),p=0$ for $v\in I^\pm_v$, respectively, so that~\eqref{eq:fastr0} admits a two-dimensional critical manifold $\mathcal{M}_0:=\{p=0, F(u,v)\}$, consisting of (at least) two branches
\begin{align}\label{eq:M0branches}
\mathcal{M}^-_0 = \{p=0, u=f^-(v), v\in I_v^-\}, \qquad \mathcal{M}^+_0 = \{p=0, u=f^+(v), v\in I_v^+)\},
\end{align}
where $F(u^\pm(v),v)=0$, and $V^-\in I_v^-$, $V^+\in I_v^+$. Assumption~\ref{assump:steadystates}\ref{assump:bistable} implies that $\mathcal{M}^\pm_0$ are of saddle type in their respective regions of definition. Using phase plane techniques, by appropriately adjusting the wave speed $c$, if $I_v^-\cap I_v^+\neq \emptyset$, then for any $v_*\in I_v^-\cap I_v^+$, there exists a locally unique speed $c_*$ and corresponding heteroclinic orbit $u_*(\xi)$ between the saddle branches $\mathcal{M}^-_0$ and $\mathcal{M}^+_0$ lying in the intersection $\mathcal{W}^{\mathrm{u}}(\mathcal{M}^-_0)\cap \mathcal{W}^{\mathrm{s}}(\mathcal{M}^+_0)$; see Figure~\ref{fig:diffusion_singular}.

We now rescale $\zeta=\delta \xi$ and consider the corresponding slow system
\begin{align}
\begin{split}\label{eq:slow}
\delta u_\zeta&= p\\
 \delta p_\zeta&= -cp-F(u,v)\\
v_\zeta&=  q\\
q_\zeta&= -\left((\bar{r}+\delta c)q+ G(u,v)\right).
\end{split}
\end{align}
Setting $\delta=0$, we obtain the corresponding reduced system
\begin{align}
\begin{split}
0&= p\\
0&= -cp-F(u,v)\\
v_\zeta&=  q\\
q_\zeta&= -\left(\bar{r}q+ G(u,v)\right).
\end{split}
\end{align}
in which the flow is restricted to the critical manifold $\mathcal{M}_0$. The reduced flow on each of the saddle branches $\mathcal{M}^-_0$ and $\mathcal{M}^+_0$ is given by the planar flows
\begin{align}
\begin{split}\label{eq:reduced}
v_\zeta&=  q\\
q_\zeta&= -\left(\bar{r}q+ G(f^\pm(v),v)\right).
\end{split}
\end{align}
By Assumption~\ref{assump:steadystates}\ref{assump:bistable}, $(U^+,V^+)=(f^+(V^+),V^+)$ are fixed points of the full system. Since 
\begin{align*}
    (f^\pm)'(v) = -F_v(f^\pm(v),v)/F_u(f^\pm(v),v),
\end{align*} 
the conditions~\eqref{eq:steadystate_conditions} imply that $G_u(f^\pm(V^\pm)(f^\pm)'(V^\pm)+G_v(f^\pm(V^\pm),V^\pm)<0$, and hence $(v,q)=(V^\pm,0)$ correspond to saddle fixed points of the reduced flows~\eqref{eq:reduced} on $\mathcal{M}_0^\pm$, respectively. We denote by $L^\mathrm{s,u}_\pm$ the stable/unstable manifolds of the fixed points $(V^\pm,0)$. In this setting, to ensure the existence of a traveling front for $0<\delta\ll 1$, we need to make the following assumption, which can be checked in a given system by examining the planar flows~\eqref{eq:reduced} for given $\bar{r}$ and nonlinearity $G(u,v)$ (see right panel of Figure~\ref{fig:diffusion_singular}).

\begin{Assumption}\label{assump:diffusion_existence}
The projection of the manifold $L^\mathrm{u}_-$ from $\mathcal{M}_0^-$ onto $\mathcal{M}_0^+$ transversely intersects $L^\mathrm{s}_+$ at $(v,q)=(v_*,q_*)$ for some $v_*\in I^-_v\cap I^+_v$.
\end{Assumption}



Considering the manifolds $L^\mathrm{s,u}_\pm$ as subsets of the critical manifold of $\mathcal{M}_0$ and taking the union of their fast (un)stable fibers, we can construct the singular two-dimensional stable and unstable manifolds of the equilibria $P^\pm$ in the full system as $\mathcal{W}^\mathrm{s,u}(P^\pm): = \mathcal{W}^\mathrm{s,u}(L^\mathrm{s,u}_\pm)\subset \mathcal{W}^\mathrm{s,u}(\mathcal{M}^\pm_0)$. 

A singular heteroclinic orbit between the equilibria $P^\pm$ can then be formed by concatenating slow/fast trajectories from the layer/reduced problems as follows. The solution leaves $P^-$ along the slow unstable manifold $L^\mathrm{u}_-\subset \mathcal{M}^-_0$, then departs $\mathcal{M}^-_0$ along a fast jump contained within $\mathcal{W}^\mathrm{u}(L^\mathrm{u}_-)$ at the critical jump value $v=v_*$. Provided $c=c_*$, by Assumption~\ref{assump:diffusion_existence} this fast jump lies in the intersection $\mathcal{W}^\mathrm{u}(L^\mathrm{u}_-)\cap \mathcal{W}^\mathrm{s}(L^\mathrm{s}_+)$, and the orbit then tracks $L^\mathrm{s}_+$ until reaching $P^+$; see Figure~\ref{fig:diffusion_singular}.  This sequence forms a singular orbit from which, for sufficiently small $\delta_0>0$, a solution to the full system~\eqref{eq:fast} for $0<\delta\leq \delta_0 $ with speed $c=c_*+\mathcal{O}(\delta)$ can be obtained using geometric singular perturbation theory.

\subsubsection{Case $(ii)$: $\delta= \bar{\delta}r, 0< \bar{\delta}\leq \bar{\delta}_0$}
We now consider $r\ll1$ and set $\delta= \bar{\delta}r$ in~\eqref{eq:fast}, obtaining
\begin{align}
\begin{split}\label{eq:fast_ii}
u_\xi&= p\\
 p_\xi&= -cp-F(u,v)\\
v_\xi&= \bar{\delta}r q\\
q_\xi&= -r\left((1+\bar{\delta}^2rc)q+\bar{\delta} G(u,v)\right).
\end{split}
\end{align}
which we similarly aim to analyze for all $0<\bar{\delta}<\bar{\delta}_0$ for any $\bar{\delta}_0>0$ where $\bar{\delta}_0$ is $\mathcal{O}(1)$ with respect to $r$. This system is now a $2$-fast-$2$-slow system with timescale separation parameter $r$, and the analysis then proceeds similarly as in the previous case for $\bar{\delta}$ bounded away from zero. 

However, when $\bar{\delta}$ is small, we employ a slightly different argument in order to obtain better estimates for the solution to the existence problem as $\bar{\delta}\to0$. First, we rescale $q = \bar{\delta}\tilde{q}$, to obtain the system
\begin{align}
\begin{split}\label{eq:fast_ii_res}
u_\xi&= p\\
 p_\xi&= -cp-F(u,v)\\
v_\xi&= \bar{\delta}^2r \tilde{q}\\
\tilde{q}_\xi&= -r\left((1+\bar{\delta}^2rc)\tilde{q}+ G(u,v)\right).
\end{split}
\end{align}
Setting $\tau=r \xi$, we obtain the reduced flow for $r=0$ on the manifolds $\mathcal{M}^\pm_0$
\begin{align}
\begin{split}\label{eq:fast_ii_reduced}
v_\tau&= \bar{\delta}^2\tilde{q}\\
\tilde{q}_\tau&= -\tilde{q}- G(f^\pm(v),v),
\end{split}
\end{align}
which can be analyzed as planar slow fast system with singular perturbation parameter $\bar{\delta}^2$. In particular, this allows us to easily determine the manifolds $L^\mathrm{u}_-$ and $L^\mathrm{s}_+$. Note that in the limit $\bar{\delta}\to0$, there exist normally attracting invariant manifolds $\mathcal{C}^\pm_0=\left\{\tilde{q}= - G(f^\pm(v),v) \right\}\subset \mathcal{M}^\pm_0$ with corresponding reduced flows
\begin{align}
\begin{split}\label{eq:fast_ii_reduced_reduced}
v_{\tilde{\tau}}&= -G(f^\pm(v),v)
\end{split}
\end{align}
where $\tilde{\tau} = \bar{\delta}^2\tau$. As argued in~\S\ref{sec:diffusion_existence_i}, Assumption~\ref{assump:steadystates} implies that the quantities $\kappa_\pm:=G_u(f^\pm(V^\pm)(f^\pm)'(V^\pm)+G_v(f^\pm(V^\pm),V^\pm)<0$, and hence the fixed points $v=V^\pm$ are repelling within $\mathcal{C}^\pm_0$. Since the manifolds $\mathcal{C}^\pm_0$ are normally attracting in the reduced flow~\eqref{eq:fast_ii_reduced}, they perturb to locally invariant manifolds $\mathcal{C}^\pm_{\bar{\delta}}$ for all sufficiently small $\bar{\delta}$, we deduce that the manifold $L^\mathrm{u}_-$ corresponds to the perturbed manifold $\mathcal{C}^-_{\bar{\delta}}$ given as a graph $\tilde{q}= - G(f^-(v),v)+\mathcal{O}(\bar{\delta}^2)$, and $L^\mathrm{s}_+$ corresponds to the perturbed stable fiber of $\mathcal{C}^+_{\bar{\delta}}$ which meets $\mathcal{C}^+_{\bar{\delta}}$ at $v=V^+$, and can be written as a graph $v=V^++\mathcal{O}(\bar{\delta}^2)$. For the manifolds $L^\mathrm{s}_+$ and $L^\mathrm{u}_-$ to intersect in their combined projection, as in Assumption~\ref{assump:diffusion_existence}, for all small $\bar{\delta}>0$, they do so at a point $v_*=V^++\mathcal{O}(\bar{\delta}^2)$ and $\tilde{q}_*=\mathcal{O}(\bar{\delta}^2)$, which therefore defines the critical jump value and speed $c=c_*$. Noting by $v^\pm(\tau)$ the solutions corresponding to $L^\mathrm{s}_+$ and $L^\mathrm{u}_-$ satisfying $v^\pm(0)=v_*$, we have that
\begin{align}
\begin{split}\label{eq:diffusion_i_existence_reduced_estimates}
    v ^-_\tau(\tau) &= -\bar{\delta}^2G(f^-(v^-(\tau)),v^-(\tau))+\mathcal{O}(\bar{\delta}^4), \quad \tau<0\\
    v^+_\tau(\tau) &= \mathcal{O}(\bar{\delta}^4), \quad \tau>0,
    \end{split}
\end{align}
and $v^+_\tau(\tau)$ decays with exponential rate $-1+\bar{\delta}^2\kappa_+$ to leading order as $\tau\to\infty$.

\begin{Remark}
    Assuming $V^->V^+$, without loss of generality, we note that to ensure that Assumption~\ref{assump:diffusion_existence} holds for all sufficiently small $\bar{\delta}$, it suffices to assume that $G(f^-(v),v)<0$ for all $v\in[V^+, V^-)$. This connects the structure of the existence problem in the intermediate regime to that in the strong advection regime; see~\S\ref{sec:existence_advection_dominant} and Assumption~\ref{assump:advection_existence}.
\end{Remark}

Using this construction in the region of small $\bar{\delta}$, we therefore obtain singular heteroclinic orbits for all $0<\bar{\delta}<\bar{\delta}_0$, which can be shown, using geometric singular perturbation techniques, to perturb to traveling front solutions with speed $c=c_*+\mathcal{O}(r)$ for all $0<r<r_0$ for some $r_0\ll1$. By choosing $\bar{r}_0>1$ in~\S\ref{sec:diffusion_existence_i} and $\bar{\delta}_0>1$ above, and possibly taking $r_0$ and/or $\delta_0$ smaller if necessary, we can combine these results with those of the previous section to obtain a slow-fast description of singular traveling fronts for any $0\leq r\leq r_0, 0<\delta\leq\delta_0$. The boundary between the weak and intermediate regimes occurs when $r=\mathcal{O}(\delta)$, or equivalently when $\nu =\mathcal{O}(\delta^{-1})$.

\subsection{Long wavelength (in)stability: leading order computation of $\lambda_{\mathrm{c},2}$}\label{sec:stability_diffusion_dominant}

Given the slightly different slow fast structures of the weak advection and intermediate regimes in~\S\ref{sec:existence_diffusion_dominant}, we similarly separate the computation of $\lambda_{\mathrm{c},2}$ in each case, depending on whether $\delta$ or $r$ is used as the primary singular perturbation parameter.

\subsubsection{Case $(i)$: $r= \bar{r}\delta, 0\leq \bar{r}\leq \bar{r}_0$}\label{sec:stabilityi_diffusion_dominant}
We consider the fast
\begin{align}\begin{split}
u_{\xi\xi}-cu_\xi +F_u(u_\mathrm{h}(\xi),v_\mathrm{h}(\xi))u +G_u(u_\mathrm{h}(\xi),v_\mathrm{h}(\xi))v &= 0\\
v_{\xi\xi}-\left(\bar{r}\delta+\delta^2c\right)v_\xi+\delta^2 F_v(u_\mathrm{h}(\xi),v_\mathrm{h}(\xi))u+\delta^2G_v(u_\mathrm{h}(\xi),v_\mathrm{h}(\xi))v&=0
\end{split}
\end{align}
and slow
\begin{align}\begin{split}
\delta^2u_{\zeta\zeta}-\delta cu_\zeta +F_u(u_\mathrm{h}(\zeta/\delta),v_\mathrm{h}(\zeta/\delta))u +G_u(u_\mathrm{h}(\zeta/\delta),v_\mathrm{h}(\zeta/\delta))v &= 0\\
v_{\zeta\zeta}-\left(\bar{r}+\delta c\right)v_\zeta+F_v(u_\mathrm{h}(\zeta/\delta),v_\mathrm{h}(\zeta/\delta))u+G_v(u_\mathrm{h}(\zeta/\delta),v_\mathrm{h}(\zeta/\delta))v&=0
\end{split}
\end{align}
formulations of the adjoint equation. In the fast field, to leading order $u_\mathrm{h}(\xi) = u_*(\xi), v_\mathrm{h}(\xi) = v_*$, and $c = c_*$, so that to leading order the fast system is given by 
\begin{align}\begin{split}
u_{\xi\xi}-c_*u_\xi +F_u(u_*(\xi),v_*)u +G_u(u_*(\xi),v_*)v &= 0\\
v_{\xi\xi}&=0
\end{split}
\end{align}
from which we deduce that  $v = \bar{v}_*= $ constant and $u$ satisfies
\begin{align}\begin{split}\label{eq:stab_fast_interim_step}
u_{\xi\xi}-c_*u_\xi +F_u(u_*(\xi),v_*)u &=- G_u(u_*(\xi),v_*)\bar{v}_*.
\end{split}
\end{align}
Taking the inner product of~\eqref{eq:stab_fast_interim_step} with the bounded solution $u_*'(\xi)e^{c_*\xi}$ of the reduced fast adjoint equation 
\begin{align}\begin{split}
u_{\xi\xi}-c_*u_\xi +F_u(u_*(\xi),v_*)u &=0,
\end{split}
\end{align}
implies that 
\begin{align}\begin{split}
0=\bar{v}_*\int_{-\infty}^\infty G_u(u_*(\xi),v_*)u_*'(\xi)\mathrm{d}\xi = \bar{v}_* \left(G(u_*^+,v_*) -G(u_*^-,v_*) \right),
\end{split}
\end{align}
By Assumption~\ref{assump:nondegeneracy}, $G(u_*^+,v_*) \neq G(u_*^-,v_*)$ so that $\bar{v}_* = 0$, and to leading order $v=0$ and $u = \alpha_* u_*'(\xi)e^{c_*\xi}$.

 For the slow system, we denote by $v^\pm(\zeta)$ the slow orbits of~\eqref{eq:reduced} corresponding to $L^\mathrm{u}_-$ and $L^\mathrm{s}_+$, respectively, and satisfying $v^\pm(0)=v_*$. At leading order we find that 
\begin{align}\begin{split}
F_u(f^\pm(v^\pm,v^\pm)u +G_u(f^\pm(v^\pm),v^\pm)v &= 0\\
v_{\zeta\zeta}-\bar{r}v_\zeta+F_v(f^\pm(v^\pm),v^\pm)u+G_v(f^\pm(v^\pm),v^\pm)v&=0
\end{split}
\end{align}
or equivalently
\begin{align}\begin{split}
v_{\zeta\zeta}-\bar{r}v_\zeta+\left(G_v(f^\pm(v^\pm),v^\pm)-F_v(f^\pm(v^\pm),v^\pm)\frac{G_u(f^\pm(v^\pm),v^\pm)}{F_u(f^\pm(v^\pm),v^\pm)}\right)v&=0,
\end{split}
\end{align}
from which we deduce that $v^{A,\pm}(\zeta) = \delta \alpha^\pm e^{\bar{r}\zeta}v^\pm_\zeta(\zeta)$. To match along the fast jump, we note that to ensure continuity of $v^A$, we require $\alpha^+=\alpha^-=\alpha$ since $v^+_\zeta(0)=q_*=v^-_\zeta(0)$. To account for the jump in $v^A_{\zeta}$ from the slow equations across the fast jump
\begin{align}
\left[(e^{\bar{r}\zeta}v^+_\zeta(\zeta))_\zeta-(e^{\bar{r}\zeta}v^-_\zeta(\zeta))_\zeta\right]_{\zeta=0}=v^+_{\zeta\zeta}(0)-v^-_{\zeta\zeta}(0)=G(u^+_*,v_*)-G(u^-_*,v_*),
\end{align}
we require a corresponding jump through the fast field
\begin{align}
\Delta v^A_\xi = -\delta^2\alpha_*\int_{-\infty}^\infty F_v(u_*(\xi),v_*)e^{c_*\xi}u_*'(\xi)\mathrm{d}\xi = \delta^2\alpha \left(v^+_{\zeta\zeta}(0)-v^-_{\zeta\zeta}(0)\right),
\end{align}
which implies
\begin{align}
\alpha = -\alpha_*\frac{\int_{-\infty}^\infty F_v(u_*(\xi),v_*)e^{c_*\xi}u_*'(\xi)\mathrm{d}\xi}{G(u^+_*,v_*)-G(u^-_*,v_*)}.
\end{align}
We now compute 
\begin{align*}
\int_{-\infty}^\infty v_\mathrm{h}'(\xi)v^A(\xi)\mathrm{d}\xi&=\int_{-\infty}^{-\frac{1}{\sqrt{\delta}}} v_\mathrm{h}'(\xi)v^A(\xi)\mathrm{d}\xi+\int_{-\frac{1}{\sqrt{\delta}}}^{\frac{1}{\sqrt{\delta}}} v_\mathrm{h}'(\xi)v^A(\xi)\mathrm{d}\xi+\int_{\frac{1}{\sqrt{\delta}}}^\infty v_\mathrm{h}'(\xi)v^A(\xi)\mathrm{d}\xi\\
&=\alpha\int_{-\infty}^{-\sqrt{\delta}} \delta e^{\bar{r}\zeta}v^-_\zeta(\zeta)^2\mathrm{d}\zeta+\alpha\int_{-\frac{1}{\sqrt{\delta}}}^{\frac{1}{\sqrt{\delta}}} \delta^2 q_*^2\mathrm{d}\xi+\alpha\int_{\sqrt{\delta}}^\infty \delta e^{\bar{r}\zeta}v^+_\zeta(\zeta)^2\mathrm{d}\zeta +\mathcal{O}(\alpha \delta^2)\\
&=\alpha\int_{-\infty}^{-\sqrt{\delta}} \delta e^{\bar{r}\zeta}v^-_\zeta(\zeta)^2\mathrm{d}\zeta+\alpha\int_{\sqrt{\delta}}^\infty \delta e^{\bar{r}\zeta}v^+_\zeta(\zeta)^2\mathrm{d}\zeta +\mathcal{O}(\alpha \delta^{3/2})
\end{align*}
and similarly
\begin{align*}
\int_{-\infty}^\infty u_\mathrm{h}'(\xi)u^A(\xi)\mathrm{d}\xi&=\alpha_*\left(\int_{-\infty}^\infty e^{c_*\xi}u_*'(\xi)^2\mathrm{d}\xi+\mathcal{O}(\delta)\right)
\end{align*}
so that to leading order, we estimate~\eqref{eq:lambda2c_expression} as
\begin{align}\label{eq:lambda2c_diffusion_dominant}
\lambda_{\mathrm{c},2} &\sim-\frac{1}{\delta}\frac{\int_{-\infty}^\infty F_v(u_*(\xi),v_*)e^{c_*\xi}u_*'(\xi)\mathrm{d}\xi}{G(u^+_*,v_*)-G(u^-_*,v_*)}\frac{\int_{-\infty}^{0} e^{\bar{r}\zeta}v^-_\zeta(\zeta)^2\mathrm{d}\zeta+ \int_{0}^\infty e^{\bar{r}\zeta}v^+_\zeta(\zeta)^2\mathrm{d}\zeta}{\int_{-\infty}^\infty e^{c_*\xi}u_*'(\xi)^2\mathrm{d}\xi}.
\end{align}
The last factor has fixed sign so that we immediately obtain the stability criterion~\eqref{eq:diffusion_criterion} in the weak advection regime in terms of the quantities~\eqref{eq:fstar_gstar}. The leading order asymptotic formula~\eqref{eq:lambda2c_diffusion_dominant} holds provided $r= \bar{r}\delta, 0\leq \bar{r}\leq \bar{r}_0, 0<\delta\leq \delta$, or in terms of $\nu$, provided $0\leq \nu \leq \tfrac{\bar{r}_0}{\delta}$. We note that the case $\nu=0$ corresponds to setting $\bar{r}=0$ in~\eqref{eq:lambda2c_diffusion_dominant}, which matches the expression obtained for the cofficient $\lambda_{\mathrm{c},2}$ in the absence of advection in~\cite[\S2]{CDLOR}.

\subsubsection{Case $(ii)$: $\delta= \bar{\delta}r, 0< \bar{\delta}\leq \bar{\delta}_0$}\label{sec:stability_intermediate}

We consider the fast
\begin{align}\begin{split}
u_{\xi\xi}-cu_\xi +F_u(u_\mathrm{h}(\xi),v_\mathrm{h}(\xi))u +G_u(u_\mathrm{h}(\xi),v_\mathrm{h}(\xi))v &= 0\\
v_{\xi\xi}-\left(r+\bar{\delta}^2r^2c\right)v_\xi+\bar{\delta}^2r^2 F_v(u_\mathrm{h}(\xi),v_\mathrm{h}(\xi))u+ \bar{\delta}^2r^2G_v(u_\mathrm{h}(\xi),v_\mathrm{h}(\xi))v&=0
\end{split}
\end{align}
and slow
\begin{align}\begin{split}
r^2u_{\tau\tau}-rcu_\tau +F_u(u_\mathrm{h}(\tau/r),v_\mathrm{h}(\tau/r))u +G_u(u_\mathrm{h}(\tau/r),v_\mathrm{h}(\tau/r))v &= 0\\
v_{\tau\tau}-\left(1+\bar{\delta}^2rc\right)v_\tau+\bar{\delta}^2 F_v(u_\mathrm{h}(\tau/r),v_\mathrm{h}(\tau/r))u+\bar{\delta}^2G_v(u_\mathrm{h}(\tau/r),v_\mathrm{h}(\tau/r))v&=0
\end{split}
\end{align}
formulations of the adjoint equation. At leading order the fast system is given by 
\begin{align}\begin{split}
u_{\xi\xi}-c_*u_\xi +F_u(u_*(\xi),v_*)u +G_u(u_*(\xi),v_*)v &= 0\\
v_{\xi\xi}&=0
\end{split}
\end{align}
from which we deduce as in~\S\ref{sec:stabilityi_diffusion_dominant} that to leading order $v=0$ and $u(\xi) = \alpha_*e^{c_*\xi}u_*'(\xi)$. For the slow system, at leading order we find that 
\begin{align}\begin{split}
F_u(f^\pm(v^\pm),v^\pm)u +G_u(f^\pm(v^\pm),v^\pm)v &= 0\\
v_{\tau\tau}-v_\tau+\bar{\delta}^2 F_v(f^\pm(v^\pm),v^\pm)u+\bar{\delta}^2G_v(f^\pm(v^\pm),v^\pm)v&=0
\end{split}
\end{align}
or equivalently
\begin{align}\begin{split}
v_{\tau\tau}-v_\tau+\bar{\delta}^2\left(G_v(f^\pm(v^\pm),v^\pm)-F_v(f^\pm(v^\pm),v^\pm)\frac{G_u(f^\pm(v^\pm),v^\pm)}{F_u(f^\pm(v^\pm),v^\pm)}\right)v&=0,
\end{split}
\end{align}
from which we deduce that $v^{A,\pm}(\tau) = r\alpha^\pm e^\tau v^\pm_\tau(\tau)$ in the slow fields. To match along the fast jump, we note that to ensure continuity of $v^A$, we require $\alpha^+=\alpha^-=\alpha$ since $v^+_\tau(0)=v^-_\tau(0)$. To account for the jump in $v^A_{\tau}$ from the slow equations across the fast jump
\begin{align}
\left[(e^\tau v^+_\tau(\tau))_\tau-(e^\tau v^-_\tau(\tau))_\tau\right]_{\tau=0}=v^+_{\tau\tau}(0)-v^-_{\tau\tau}(0)=\bar{\delta}^2\left(G(u^+_*,v_*)-G(u^-_*,v_*)\right),
\end{align}
we require a corresponding jump through the fast field
\begin{align}
\Delta v^A_\xi = -\bar{\delta}^2r^2\alpha_*\int_{-\infty}^\infty F_v(u_*(\xi),v_*)e^{c_*\xi}u_*'(\xi)\mathrm{d}\xi = r^2\alpha \left(v^+_{\tau\tau}(0)-v^-_{\tau\tau}(0)\right),
\end{align}
which implies
\begin{align}
\alpha = -\alpha_*\frac{\int_{-\infty}^\infty F_v(u_*(\xi),v_*)e^{c_*\xi}u_*'(\xi)\mathrm{d}\xi}{G(u^+_*,v_*)-G(u^-_*,v_*)}.
\end{align}
Proceeding as in~\S\ref{sec:stabilityi_diffusion_dominant}, we estimate
\begin{align*}
\int_{-\infty}^\infty v_\mathrm{h}'(\xi)v^A(\xi)\mathrm{d}\xi&=\alpha\int_{-\infty}^{-\sqrt{r}} r e^{\tau}v^-_\tau(\tau)^2\mathrm{d}\tau+\alpha\int_{\sqrt{r}}^\infty r e^{\tau}v^+_\tau(\tau)^2\mathrm{d}\tau +\mathcal{O}(\alpha r^{3/2})\\
\int_{-\infty}^\infty u_\mathrm{h}'(\xi)u^A(\xi)\mathrm{d}\xi&=\alpha_*\left(\int_{-\infty}^\infty e^{c_*\xi}u_*'(\xi)^2\mathrm{d}\xi+\mathcal{O}(r)\right)
\end{align*}

We now estimate~\eqref{eq:lambda2c_expression} to leading order as
 \begin{align*}
\lambda_{\mathrm{c},2} &\sim -1+\frac{\alpha r}{\alpha_*\delta^2}\frac{\int_{-\infty}^{0} e^\tau v^-_\tau(\tau)^2\mathrm{d}\tau+ \int_{0}^\infty e^\tau v^+_\tau(\tau)^2\mathrm{d}\tau+\mathcal{O}(\sqrt{r})}{\int_{-\infty}^\infty e^{c_*\xi}u_*'(\xi)^2\mathrm{d}\xi}\\
&=-1-\frac{1}{\bar{\delta}^2r}\frac{\int_{-\infty}^\infty F_v(u_*(\xi),v_*)e^{c_*\xi}u_*'(\xi)\mathrm{d}\xi}{G(u^+_*,v_*)-G(u^-_*,v_*)}\frac{\int_{-\infty}^{0} e^\tau v^-_\tau(\tau)^2\mathrm{d}\tau+ \int_{0}^\infty e^\tau v^+_\tau(\tau)^2\mathrm{d}\tau+\mathcal{O}\left( \sqrt{r}\right)}{\int_{-\infty}^\infty e^{c\xi}u_f'(\xi)^2\mathrm{d}\xi}.
\end{align*}
When $\bar{\delta}$ is small, we can evaluate further using~\eqref{eq:diffusion_i_existence_reduced_estimates} so that to leading order
 \begin{align*}\label{eq:lambda2c_intermediate}
\lambda_{\mathrm{c},2} &=-1-\frac{\bar{\delta}^2}{r}\frac{\int_{-\infty}^\infty F_v(u_*(\xi),v_*)e^{c_*\xi}u_*'(\xi)\mathrm{d}\xi}{G(u^+_*,v_*)-G(u^-_*,v_*)}\frac{\int_{-\infty}^{0} e^\tau
G(f^-(v^-(\tau)), v^-(\tau))^2\mathrm{d}\tau}{\int_{-\infty}^\infty e^{c_*\xi}u_*'(\xi)^2\mathrm{d}\xi}+\mathcal{O}\left( \frac{\bar{\delta}^2}{\sqrt{r}}, \frac{\bar{\delta}^4}{r}\right).
\end{align*}
From this, we obtain the stability criterion~\eqref{eq:intermediate_criterion} in the intermediate regime. The leading order expression~\eqref{eq:intermediate_criterion} is valid provided $\delta= \bar{\delta}r, 0< \bar{\delta}\leq \bar{\delta}_0, 0\leq r\leq r_0$, or in terms of $\nu$, provided $\tfrac{1}{\bar{\delta}_0\delta}\leq \nu \leq \tfrac{r_0}{\delta^2}$. We see that the first (constant) term $-1$ dominates when $\bar{\delta}^2\ll r$, while the second term dominates when $r\ll\bar{\delta}^2$. If the coefficient of the latter is positive, then the expression~\eqref{eq:lambda2c_intermediate} changes sign when
 \begin{align*}
r&=-\bar{\delta}^2\frac{\int_{-\infty}^\infty F_v(u_*(\xi),v_*)e^{c_*\xi}u_*'(\xi)\mathrm{d}\xi}{G(u^+_*,v_*)-G(u^-_*,v_*)}\frac{\int_{-\infty}^{0} e^\tau
G(f^-(v^-(\tau)), v^-(\tau)^2\mathrm{d}\tau}{\int_{-\infty}^\infty e^{c_*\xi}u_*'(\xi)^2\mathrm{d}\xi}+\mathcal{O}\left(\bar{\delta}^4\right),
\end{align*}
or equivalently when
 \begin{align*}
\nu=M^{1/3}\delta^{-4/3}\left(1+\mathcal{O}(\delta^{2/3})\right)
\end{align*}
where
\begin{align}
M:=\frac{-\int_{-\infty}^\infty F_v(u_*(\xi),v_*)e^{c_*\xi}u_*'(\xi)\mathrm{d}\xi}{G(u^+_*,v_*)-G(u^-_*,v_*)}\frac{\int_{-\infty}^{0} e^\tau
G(f^-(v^-(\tau)), v^-(\tau))\mathrm{d}\tau}{\int_{-\infty}^\infty e^{c_*\xi}u_*'(\xi)^2\mathrm{d}\xi}.
\end{align}
We emphasize that the sign of $M$ is determined by the same quantities $F_*, G_*$ defined in~\eqref{eq:fstar_gstar} which determine the sign of $\lambda_{\mathrm{c},2}$ in the weak advection regime.

\section{Strong advection regime}\label{sec:advection_dominant}

We now determine the behavior of the system in the strong advection regime for $r\geq r_0$, where $r_0$ is the small constant from~\S\ref{sec:diffusion_dominant}, or equivalently, for $\nu\geq \tfrac{r_0}{\delta^2}$.

\subsection{Slow/fast structure of the traveling front}\label{sec:existence_advection_dominant}

We rescale $q=\frac{\delta}{r}\tilde{q}$ and obtain
\begin{align}
\begin{split}
u_\xi&= p\\
p_\xi&= -cp -F(u,v)\\
v_\xi&= \eps \tilde{q}\\
 \tilde{q}_\xi&= r\left(-\left(1+\eps c\right)\tilde{q} -G(u,v)\right),
\end{split}
\end{align}
where we recall $\eps := \nu^{-1}=\frac{\delta^2}{r}$. We now consider this equation in the regime $0<\eps \ll1$ and $r\geq r_0$. We use $\eps$ as the timescale separation parameter, so that $v$ is a slow variable, and $(u,p,\tilde{q})$ are fast variables. We note that if $r$ is large, then $\tilde{q}$ evolves on a third (``super" fast) timescale. Setting $\eps=0$ defines the layer problem
\begin{align}
\begin{split}\label{eq:strongad_layer}
u_\xi&= p\\
p_\xi&= -cp -F(u,v)\\
 \tilde{q}_\xi&= -r\left(\tilde{q} +G(u,v)\right),
\end{split}
\end{align}
and the critical manifold
\begin{align}
\tilde{\mathcal{M}}_0:=\{p=0, F(u,v)=0, \tilde{q}= -G(u,v)\},
\end{align}
which admits at least two normally hyperbolic branches $\tilde{\mathcal{M}}_0^\pm=\left\{p=0, u=f^\pm(v), \tilde{q}= -G(f^\pm(v),v)\right\}$. The reduced flow on these manifolds is given by
\begin{align}
\begin{split}\label{eq:reduced_int}
v_\eta&= - G(f^\pm(v),v)
\end{split}
\end{align}
with respect to $\eta = \eps \xi$. By Assumption~\ref{assump:steadystates}, $G_u(f^\pm(V^\pm)(f^\pm)'(V^\pm)+G_v(f^\pm(V^\pm),V^\pm)<0$, so that the equilibria $v=V^\pm$ are repelling on $\tilde{\mathcal{M}}^\pm_0$. We assume the following.
\begin{Assumption}\label{assump:advection_existence}
Without loss of generality, we assume $V^->V^+$ and that $[V^+,V^-]\subseteq I_v^-$. Furthermore, $G(f^-(v),v)<0$ for $v\in [V^+,V^-)$.
\end{Assumption}

\begin{figure}[t]
\centering
\includegraphics[width=0.4\linewidth]{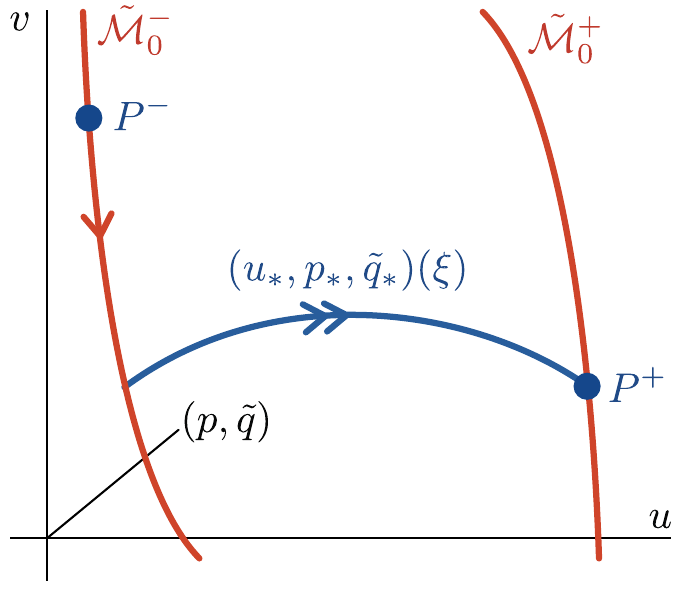}
\caption{Schematic of the slow/fast construction of the singular traveling front in the strong advection regime.  }
\label{fig:advection_singular}
\end{figure}
In this case, the unstable manifold of the equilibrium $v=V^-$ in the reduced flow~\eqref{eq:reduced_int} includes the segment $\tilde{\mathcal{M}}^-_0\cap\{v^+\leq v<v^-\}$, consisting of a single solution; we denote by $v=v^-(\eta)$ the corresponding solution satisfying $v^-(0)=v^+$. On the other hand, the stable manifold of $v=V^+$ in the reduced flow~\eqref{eq:reduced_int} on $\tilde{\mathcal{M}}^+_0$ is trivial; see Figure~\ref{fig:advection_singular}. Thus, to form a singular heteroclinic orbit between $P^\pm$ in the full system, it is necessary for the fast jump to occur in the subspace $v=v^+$, which can be arranged by taking $v_*=v^+$, and choosing $c_*$ appropriately to obtain a fast heteroclinic orbit $(u_*,p_*)(\xi)$ in the $(u,p)$-subsystem of the layer problem
\begin{align}
\begin{split}
u_\xi&= p\\
 p_\xi&= -cp-F(u,v_*).
\end{split}
\end{align}
The associated $\tilde{q}$-profile can then be obtained by integrating
\begin{align}
\begin{split}\label{eq:frontsp}
\tilde{q}_*(\xi) &= -r\int_{-\infty}^\xi e^{-r(\xi-\tilde{\xi})}G(u_*(\tilde{\xi}),v_*)\mathrm{d}\tilde{\xi}.
\end{split}
\end{align}
We note that the additional fast direction adds a uniformly stable hyperbolic direction to each of the critical manifolds $\bar{\mathcal{M}}^0_0$, $\bar{\mathcal{M}}^+_0$, and hence just adds $+1$ to the dimension of their respective stable manifolds. 

To make this construction uniform in the limit $r\to\infty$ (so that we obtain a uniform slow/fast description of the traveling front for all $r\geq r_0$ and $\eps\leq \eps_0$ for some $\eps_0$ independent of $r$), we note that in the limit $r\to\infty$, the layer problem~\eqref{eq:strongad_layer} is itself a slow/fast system with timescale separation parameter $1/r$. The existence of a fast heteroclinic orbit in the layer problem can be obtained (uniformly) for all sufficiently large $r$ using geometric singular perturbation theory. Hence the slow fast structure is well defined for any sufficiently small $\eps$, any $r\geq r_0$, and any sufficiently small $\delta$. However, due to the relation between $r,\delta,\eps$, this region includes all $r\geq  r_0$, $\delta\leq \delta_0$, and $\eps<\eps_0$ for some $\eps_0,\delta_0\ll1$, where the constant $\delta_0$ from~\S\ref{sec:diffusion_dominant} may be taken smaller, if necessary.


\subsection{Long wavelength stability: leading-order computation of $\lambda_{\mathrm{c},2}$}\label{sec:stability_advection_dominant}

We now consider the fronts from~\S\ref{sec:existence_advection_dominant}, using $\eps$ as the timescale separation parameter. We consider the fast
\begin{align}\begin{split}\label{eq:intermediate_fast}
u_{\xi\xi}-cu_\xi +F_u(u_\mathrm{h}(\xi),v_\mathrm{h}(\xi))u +G_u(u_\mathrm{h}(\xi),v_\mathrm{h}(\xi))v &= 0\\
\frac{1}{r}v_{\xi\xi}-\left(1+\eps c\right)v_\xi+\eps  F_v(u_\mathrm{h}(\xi),v_\mathrm{h}(\xi))u+ \eps G_v(u_\mathrm{h}(\xi),v_\mathrm{h}(\xi))v&=0
\end{split}
\end{align}
and slow
\begin{align}\begin{split}
\eps^2u_{\eta\eta}-\eps cu_\eta +F_u(u_\mathrm{h}(\eta/\eps),v_\mathrm{h}(\eta/\eps))u +G_u(u_\mathrm{h}(\eta/\eps),v_\mathrm{h}(\eta/\eps))v &= 0\\
\frac{\eps }{r}v_{\eta\eta}-\left(1+\eps c\right)v_\eta+F_v(u_\mathrm{h}(\eta/\eps),v_\mathrm{h}(\eta/\eps))u+G_v(u_\mathrm{h}(\eta/\eps),v_\mathrm{h}(\eta/\eps))v&=0
\end{split}
\end{align}
formulations of the adjoint equation, where $\eta = \eps\xi$. Proceeding as in~\S\ref{sec:stabilityi_diffusion_dominant}, in the fast field, to leading order $u_\mathrm{h}(\xi) = u_*(\xi), v_\mathrm{h}(\xi) = v_*$, and $c = c_*$, so that  $v = \bar{v}_*= $ constant and $u$ satisfies
\begin{align}\begin{split}
u_{\xi\xi}-c_*u_\xi +F_u(u_*(\xi),v_*)u &=- G_u(u_*(\xi),v_*)\bar{v}_*,
\end{split}
\end{align}
which implies that 
\begin{align}\begin{split}
0=\bar{v}_*\int_{-\infty}^\infty G_u(u_*(\xi),v_*)u_*'(\xi)\mathrm{d}\xi = \bar{v}_* \left(G(u_*^+,v_*) -G(u_*^-,v_*) \right),
\end{split}
\end{align}
By Assumption~\ref{assump:nondegeneracy}, $G(u_*^+,v_*) \neq G(u_*^-,v_*)$ so that $\bar{v}_* = 0$ to leading order, and $u = \alpha_* u_*'(\xi)e^{c_*\xi}$. In the slow field (recall there is only a slow field along $\mathcal{M}^-_0$), to leading order, we have 
\begin{align}\begin{split}
F_u(f^-(v^-),v^-)u +G_u(f^-(v^-),v^-)v &= 0\\
-v_\eta+F_v(f^-(v^-),v^-)u+G_v(f^-(v^-),v^-)v&=0
\end{split}
\end{align}
or equivalently
\begin{align}\begin{split}
-v_\eta+\left(G_v(f^-(v^-),v^-) -F_v(f^-(v^-),v^-)\frac{G_u(f^-(v^-),v^-)}{F_u(f^-(v^-),v^-)}  \right)v&=0,
\end{split}
\end{align}
so that $v(\eta) = \alpha^- v^-(\eta)^{-1}$ in the slow field to leading order, from which we deduce that $\alpha^-=0$ and hence $v(\eta)=\mathcal{O}(\eps)$ in the slow field. Thus, we obtain (at leading order)
\begin{align*}
\lambda_{\mathrm{c},2} &= -\frac{\int_{-\infty}^\infty \alpha_*u_*'(\xi)^2e^{c_*\xi}\mathrm{d}\xi}{\int_{-\infty}^\infty \alpha_*u_*'(\xi)^2e^{c_*\xi}\mathrm{d}\xi}+\mathcal{O}\left(\frac{\eps}{r}\right)=-1+\mathcal{O}\left(\frac{\eps}{r}\right),
\end{align*}
which directly implies~\eqref{eq:advection_criterion}.

\section{Application to a dryland ecosystem  model}\label{sec:dryland_example}
We apply the results of~\S\ref{sec:diffusion_dominant}--\ref{sec:advection_dominant} to the following modified Klausmeier model~\cite{klausmeier1999regular} in the specific form introduced in~\cite{BCD}
\begin{align}
\begin{split}\label{eq:modifiedKlausmeier}
U_t &=\Delta U-\mu_1U+U^2V\left(1-\mu_2U\right)\\
V_t &= \frac{1}{\delta^2} \Delta V +\mu_3-V-U^2V+\nu V_x,
\end{split}
\end{align}
corresponding to~\eqref{eq:grda} with $\boldsymbol{\mu}=(\mu_1, \mu_2, \mu_3)$ and 
\begin{align}
\begin{split}\label{eq:grda_mk}
F(U,V;\boldsymbol{\mu})&= -\mu_1U+U^2V\left(1-\mu_2 U\right),\\
G(U,V;\boldsymbol{\mu})&= \mu_3-V-U^2V.
\end{split}
\end{align}
 In the context of dryland ecosystem $U$ represents biomass of a species of vegetation, while $V$ represents a limiting resource such as water. The model~\eqref{eq:modifiedKlausmeier} also corresponds to that considered in~\cite{eigentler2021species} in the case of a single species. The system parameters $\mu_1, \mu_2, \mu_3$ are positive and represent mortality, inverse of soil carrying capacity, and rainfall, respectively. The small parameter $0<\delta\ll1$ representing the ratio of diffusion coefficients reflects the fact that water diffuses more quickly than vegetation, while the advection term $\nu V_x$ models the downhill flow of water on a (constantly) sloped terrain, whose slope is oriented in the $x$-direction, with the coefficient $\nu$ describing the grade of the slope. The primary differences between the model~\eqref{eq:modifiedKlausmeier} and Klausmeier's original model~\cite{klausmeier1999regular} are the inclusion of the large diffusion term in the water ($V$) equation, and the additional factor $(1-\mu_2 U)$ in the nonlinearity in the $U$ equation, representing the carrying capacity of soil.

\begin{Remark}
In~\cite{BCD}, the parameter $\nu$ was assumed large, which reflects the comparatively fast timescale on which water $V$ flows downhill (modeled via the advection term in the $u$ equation) compared to the rate at which the vegetation $U$ diffuses. It was shown in~\cite{BCD} in the absence of the diffusion term in the $V$ equation, i.e. in the limit $\delta\to\infty$ that stable planar vegetation fronts and stripes can form, aligned in the direction transverse to the slope of the terrain. Ecologically, one expects, however, that water does diffuse due to soil-water transport, and at a rate faster than that of vegetation~\cite{gilad2004ecosystem}, so that in fact one should include a large diffusion term in the $V$ equation, as in~\eqref{eq:modifiedKlausmeier}.
\end{Remark}

This system was analyzed in the absence of advection, i.e. $\nu=0$ in~\cite[\S4.1]{CDLOR}, where it was shown that bistable interfaces are \emph{always} unstable to long wavelength perturbations. Using the results of~\S\ref{sec:diffusion_dominant}--\ref{sec:advection_dominant}, we show below that this same instability is present throughout the weak advection regime, but that the presence of sufficiently large ($\nu\gg1$) advection can stabilize the fronts to long wavelength perturbations. 

\begin{figure}
\centering
\includegraphics[width=0.4\linewidth]{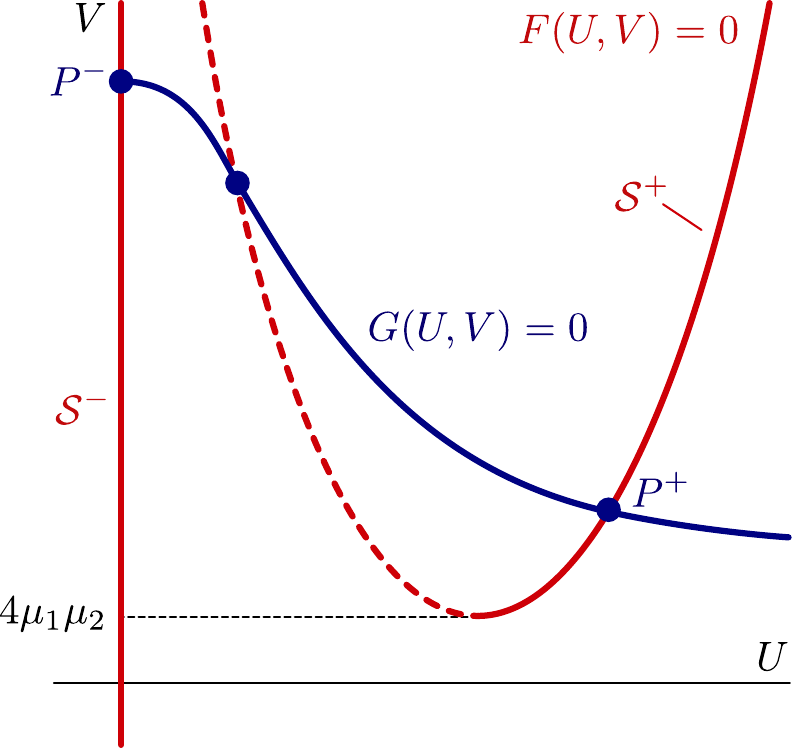}
\caption{  Sketch of the nullclines of~\eqref{eq:grda_mk} in the case~\eqref{eq:mk_steadystate_stable}, depicting the branches $\mathcal{S}^\pm$ and equilibria $P^\pm$.}
\label{fig:klaus_nullclines}
\end{figure}

In~\eqref{eq:modifiedKlausmeier}, the set $\mathcal{S}=\{F(U,V;\boldsymbol{\mu})=0\}$ is comprised of three branches
\begin{align}
U=0, \qquad U=U^\pm_F(V):= \frac{1\pm\sqrt{1-\frac{4\mu_1 \mu_2}{ V}}}{2\mu_2}
\end{align}
where the latter two branches are defined in the region $V\geq 4\mu_1\mu_2$. To find steady states, we solve
\begin{align}
    F(U,V;\boldsymbol{\mu})=G(U,V;\boldsymbol{\mu})=0,
\end{align}
and deduce that~\eqref{eq:modifiedKlausmeier} admits a stable steady state $(U^-, V^-):= (0,\mu_3)$ corresponding to zero vegetation (the `desert' state). When
\begin{align}
\frac{\mu_3}{\mu_1}> 2\left( \mu_2+\sqrt{1+\mu_2^2}\right)
\end{align}
there are two additional steady states $(U,V)=(U_{1,2},V_{1,2})$ representing uniform vegetation where
\begin{align*}
U_{1,2}&= \frac{\mu_3\pm \sqrt{\mu_3^2-4\mu_1(\mu_1+\mu_2\mu_3)}}{2(\mu_1+\mu_2\mu_3)}, \qquad V_{1,2}= \mu_3 -\frac{\mu_1 U_{1,2}}{1-\mu_2 U_{1,2}},
\end{align*}
Furthermore, when 
\begin{align}\label{eq:mk_steadystate_stable}
\frac{\mu_3}{\mu_1}> 4\mu_2+\frac{1}{\mu_2}
\end{align}
the state $(U^+,V^+):=(U_2,V_2)$ is PDE stable, i.e. the condition~\ref{assump:steadystates}\ref{assump:bistable} is satisfied, while the remaining steady state $(U_1,V_1)$ is always unstable~\cite{BCD}. Within the set $\mathcal{S}$, the branches
\begin{align}
    \mathcal{S}^-:=\left\{U=f^-(V)\equiv 0\right\},\qquad \mathcal{S}^+:= \left\{U=f^+(V)\equiv U^+_F(V), \right\}
\end{align}
contain the steady states of interest: the state $(U^-, V^-)$ always lies on the branch $\mathcal{S}^-$, and when~\eqref{eq:mk_steadystate_stable} is satisfied, the state $(U^+, V^+)$ resides on $\mathcal{S}^+$; see Figure~\ref{fig:klaus_nullclines}. When considering bistable traveling fronts, we are therefore interested in fronts between the steady states $(U^-,V^-)$ and $(U^+,V^+)$ in the parameter regime~\eqref{eq:mk_steadystate_stable} which traverse a singular heteroclinic orbit of the fast layer equation
\begin{align}
\begin{split}\label{eq:modifiedKlausmeier_layer}
0&=U_{\xi\xi}+cU_\xi-\mu_1U+U^2V\left(1-\mu_2U\right)
\end{split}
\end{align}
between the branches $\mathcal{S}^-$ and $\mathcal{S}^+$. We note that the cubic nonlinearity allows for the explicit construction of heteroclinic orbits in~\eqref{eq:modifiedKlausmeier_layer} between $\mathcal{S}^-$ and $\mathcal{S}^+$ for each $v_*>4\mu_1\mu_2$; see e.g.~\cite[\S2]{BCD} or~\cite[\S4.1]{CDLOR}. 

Using geometric singular perturbation techniques as in~\cite{BCD}, and separately considering the scalings in~\S\ref{sec:diffusion_dominant}-\ref{sec:advection_dominant}, it is possible to rigorous construct bistable traveling fronts between $(U^\pm, V^\pm)$ in~\eqref{eq:modifiedKlausmeier} for $0<\delta\ll 1$. Since we are primarily interested in the effect of advection on the stability of such interfaces, we do not carry out a detailed existence analysis here, but point to other works which rigorously construct traveling fronts and other traveling or stationary waves in various parameter regimes in the same equation~\eqref{eq:grda_mk}: in particular, the condition  (Assumption~\ref{assump:advection_existence}) necessary for existence in the strong advection regime  is verified in~\cite{BCD}, while in the weak advection and intermediate regimes, Assumption~\ref{assump:diffusion_existence} can be verified using similar techniques as in~\cite{BCDL}. Therefore, in the remainder we assume the existence of a family of traveling fronts $\phi_\mathrm{h}(\xi;\nu, \delta) = (u_\mathrm{h},v_\mathrm{h})(\xi;\nu, \delta)$ for $0\leq \nu <\infty$ and sufficiently small $\delta\ll1$.

\begin{figure}
\hspace{.05\textwidth}
\begin{subfigure}{.35 \textwidth}
\centering
\includegraphics[width=1\linewidth]{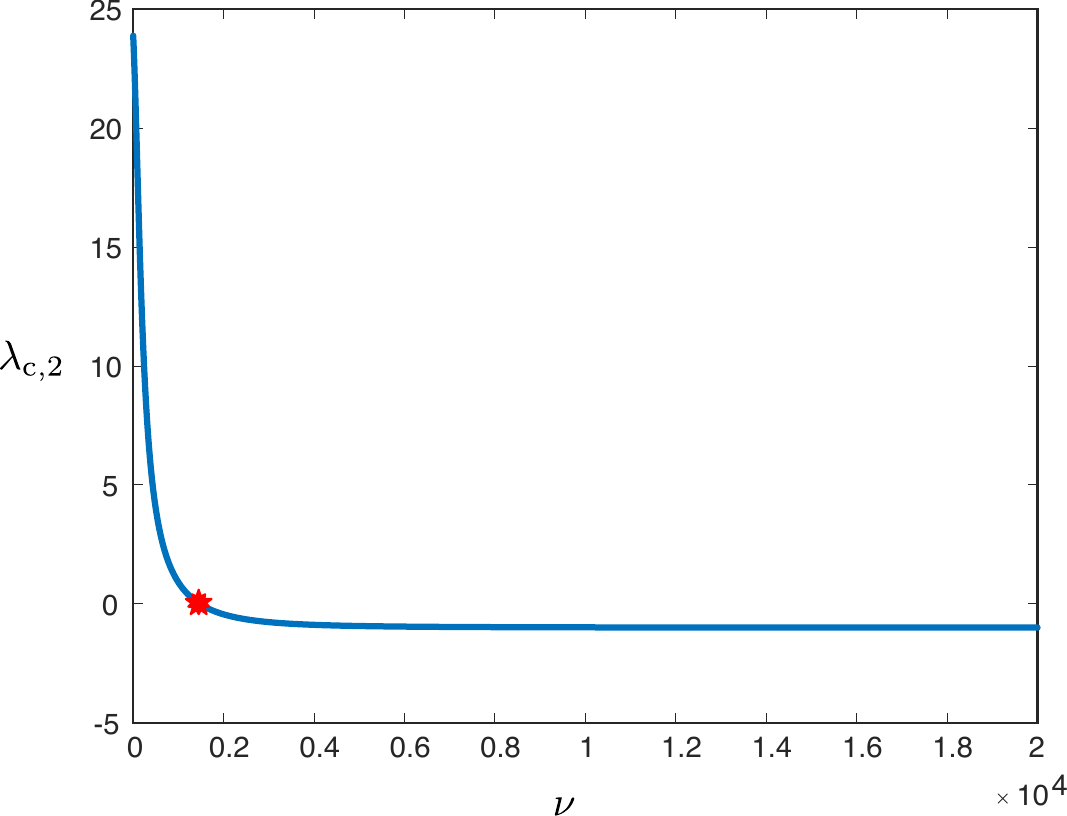}
\caption{ }
\end{subfigure}
\hspace{.1\textwidth}
\begin{subfigure}{.35 \textwidth}
\centering
\includegraphics[width=1\linewidth]{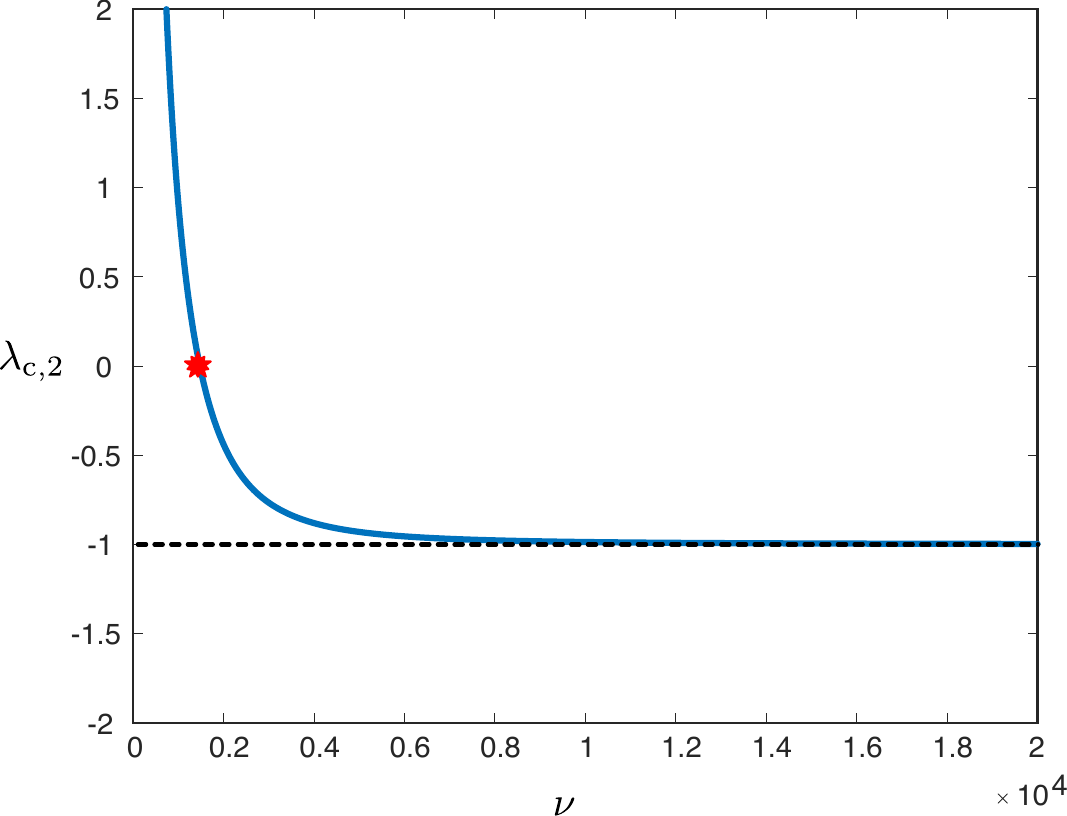}
\caption{ }
\end{subfigure}
\caption{(a) Shown is a continuation branch of traveling front solutions of~\eqref{eq:modifiedKlausmeier} over a range of $\nu \in[0, 2\cdot10^4]$ for fixed $\delta=0.001$ and $\boldsymbol{\mu}=(0.1,0.1,2.0)$. The coefficient $\lambda_{\mathrm{c},2}$ is plotted versus $\nu$. We observe that $\lambda_{\mathrm{c},2}$ changes sign as $\nu$ increases at $\nu\approx 1472$. (b) Vertical zoom of the same plot as in (a), showing the horizontal asymptote $\lambda_{\mathrm{c},2}=-1$. }
\label{fig:mk_lambda2c_cont}
\end{figure}

\begin{figure}
\centering
\includegraphics[width=0.45\linewidth]{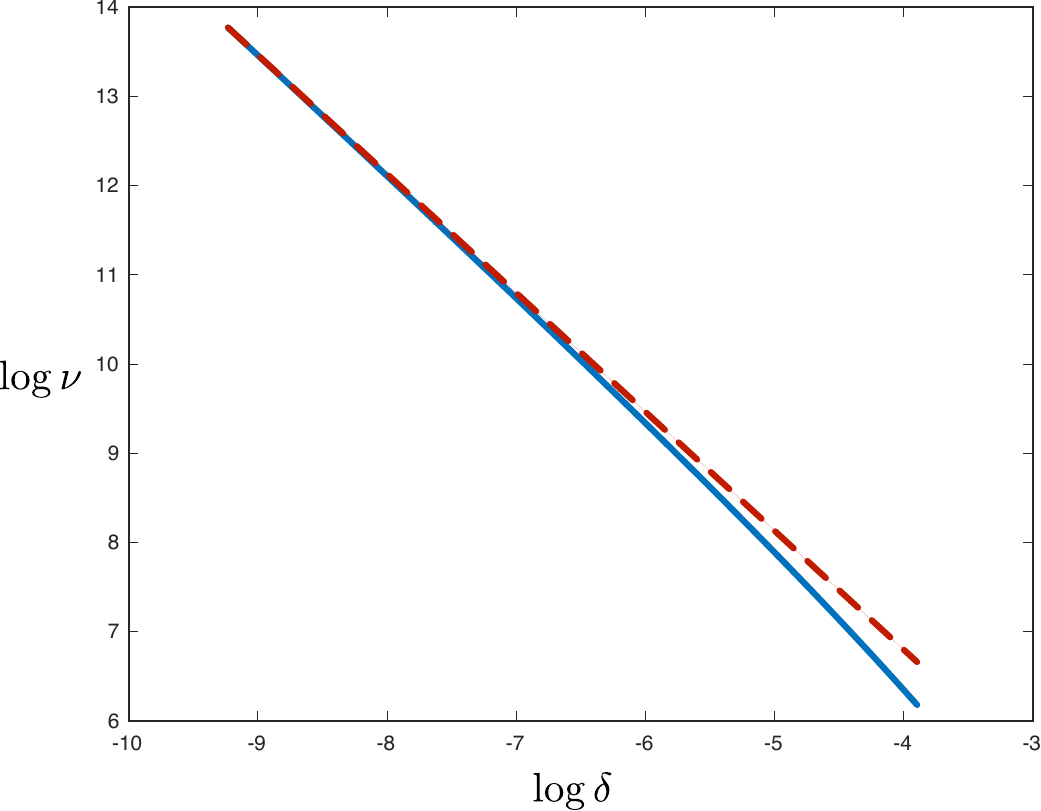}
\caption{  Results obtained by numerical continuation of the equation $\lambda_{\mathrm{c},2}=0$ for values of $\delta\in [0.0001, 0.02]$ and system parameters $\boldsymbol{\mu}=(0.1,0.1,2.0)$: a log-log plot of $\nu$ versus $\delta$ (blue) is shown alongside a line of slope $-4/3$ (dashed red).}
\label{fig:mk_zerocont}
\end{figure}

Figure~\ref{fig:mk_lambda2c_cont} depicts the results of numerically continuing the coefficient $\lambda_{\mathrm{c},2}$ using the formula~\eqref{eq:lambda2c_expression} along such a family of fronts for fixed $\delta=0.001$ and a wide range of $\nu$ values. We now examine the behavior of the coefficient $\lambda_{\mathrm{c},2}$ in the context of the asymptotic results of~\S\ref{sec:diffusion_dominant}--\ref{sec:advection_dominant}. We evaluate
      \begin{align*}
        F_v(u_*(\xi),v_*)e^{c_*\xi}&= u_*(\xi)^2\left(1-\frac{u_*(\xi)}{\mu_2}\right)>0\\
       G(u^+_*,v_*)-G(u^-_*,v_*)&=-v_*(u_*^+)^2<0
    \end{align*}
from which immediately deduce that $F_*>0$ and $G_*<0$. In particular, in the strong advection regime, this determines the sign of $\lambda_{\mathrm{c},2}$ to be positive, and thus bistable planar fronts are always unstable to long wavelength perturbations in this regime. This agrees with the results of~\cite[\S4.1]{CDLOR} in which the same equation was analyzed in the absence of advection ($\nu=0$). In the strong advection regime, based on the results of~\S\ref{sec:advection_dominant}, we find that the coefficient $\lambda_{\mathrm{c},2}$ eventually changes sign as $\nu$ increases, so that suitably large advection therefore stabilizes the fronts to long wavelength perturbations. Comparing with Figure~\ref{fig:mk_lambda2c_cont}, which plots the coefficient $\lambda_{\mathrm{c},2}$ as a function of $\nu$, and we observe that its sign indeed changes upon increasing $\nu$ at $\nu\approx 1472$. We recall that the asymptotic results of~\S\ref{sec:stability_intermediate} predict the change of sign to occur when $\nu\sim \delta^{-4/3}$. In Figure~\ref{fig:mk_zerocont}, we plot the results of continuing the equation $\lambda_{\mathrm{c},2}=0$ in the parameters $\nu, \delta$. A log-log plot of $\nu$ versus $\delta$ demonstrates good agreement with the predicted $-4/3$ exponent. For larger values of $\nu$, we see in Figure~\ref{fig:mk_lambda2c_cont} that $\lambda_{\mathrm{c},2}$ remains negative and approaches $-1$ as $\nu \to \infty$, which further agrees with the asymptotic predictions in the strong advection regime in~\S\ref{sec:stability_advection_dominant}. 

\section{Discussion}\label{sec:discussion}

In this work, we examined the effect of advection on the stability of planar fronts in two spatial dimensions, motivated by observations of planar interfaces between desert and vegetated states in water limited ecosystems. In particular, we demonstrated a stabilizing mechanism of advection on a critical eigenvalue associated with long wavelength perturbations. In the context of application to dryalnd ecosystem models, this matches the common observation that interfaces between vegetation and bare soil present in vegetation patterns tend to align perpendicular to sloped terrain, on sufficiently steep slopes. The stability criteria posed in~\S\ref{sec:summary} are model-independent, which allows for straightforward application of the results to a wide class of reaction-diffusion-advection models, under modest assumptions.


\begin{figure}
\hspace{.01\textwidth}
\begin{subfigure}{.22 \textwidth}
\centering
\includegraphics[width=1\linewidth]{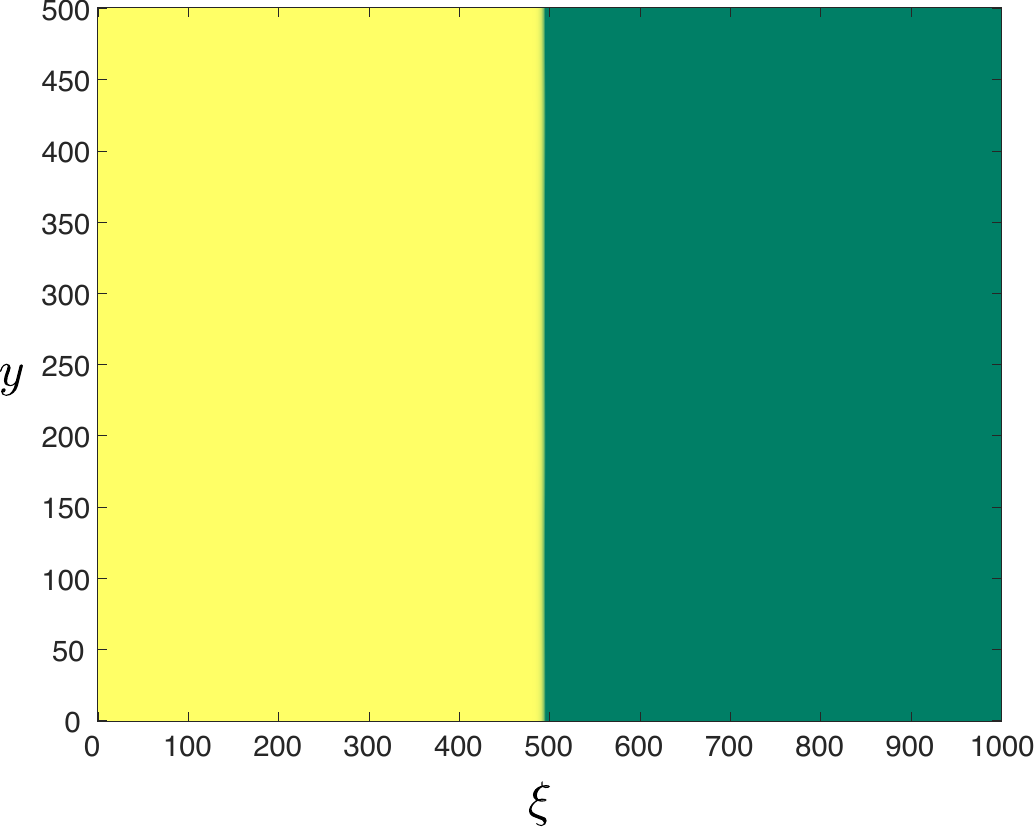}
\end{subfigure}
\hspace{.01\textwidth}
\begin{subfigure}{.22 \textwidth}
\centering
\includegraphics[width=1\linewidth]{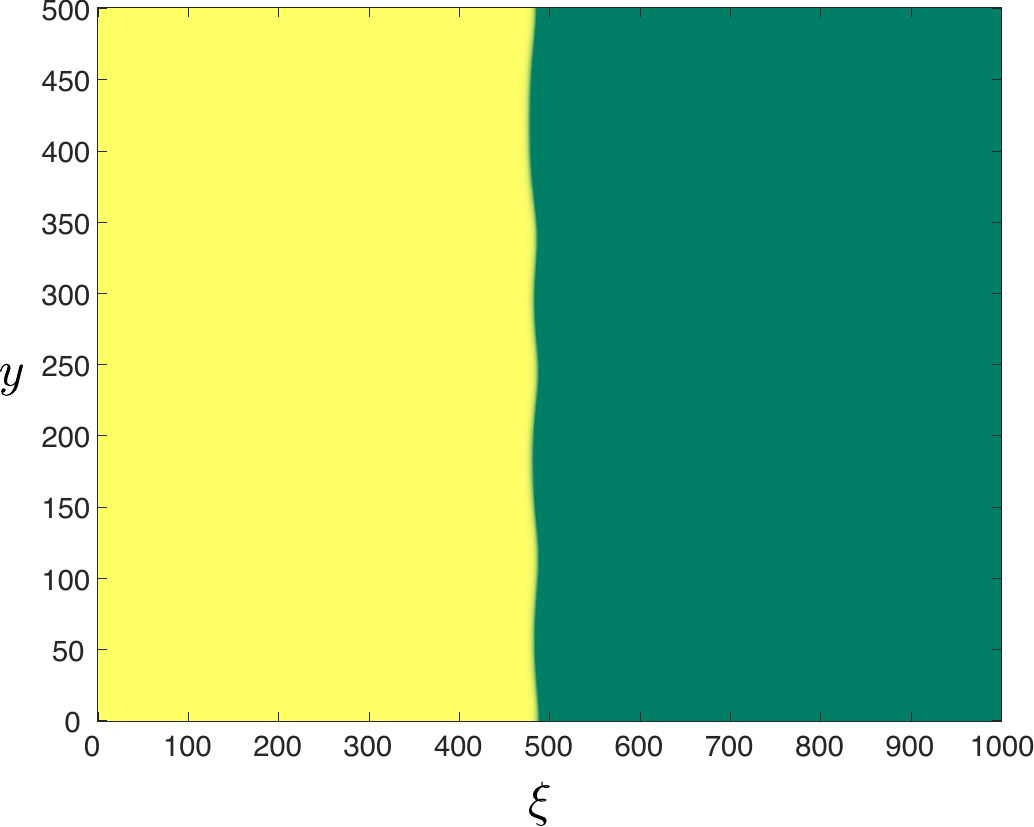}
\end{subfigure}
\hspace{.01\textwidth}
\begin{subfigure}{.22 \textwidth}
\centering
\includegraphics[width=1\linewidth]{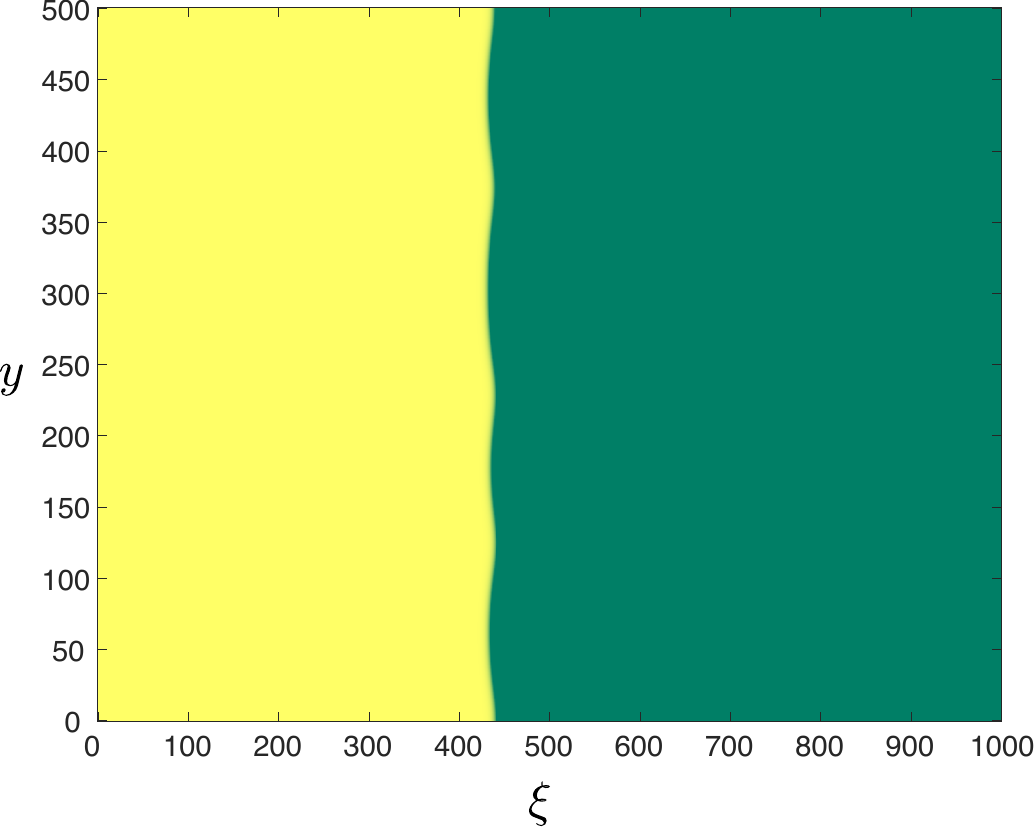}
\end{subfigure}
\hspace{.01\textwidth}
\begin{subfigure}{.22 \textwidth}
\centering
\includegraphics[width=1\linewidth]{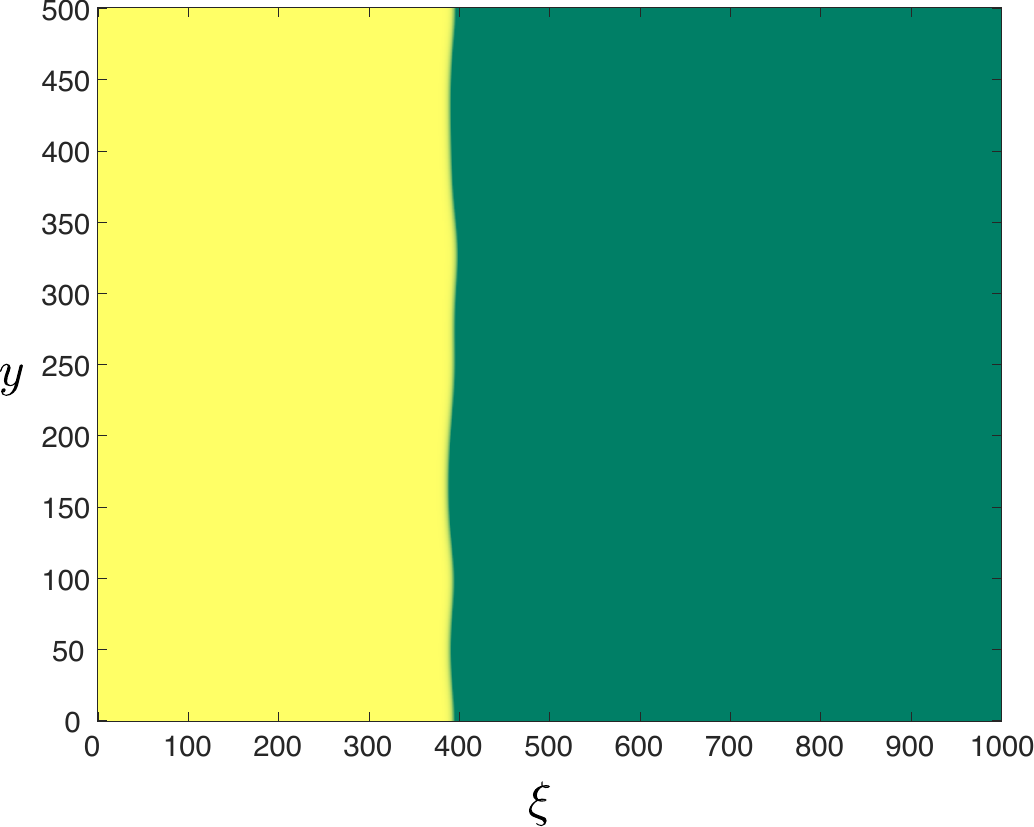}
\end{subfigure}
\caption{Shown are the results of direct numerical simulation in~\eqref{eq:grda_mk} for $\delta=0.01$, $\boldsymbol{\mu}=(0.1,0.1,2.0)$ and $\nu=1200$, which is below the critical stability boundary at which $\lambda_{\mathrm{c},2}$ changes sign (estimated at $\nu\approx 1472$). The simulation was initialized with a perturbation of the corresponding traveling front $\phi_\mathrm{h}(\xi;\nu,\delta)$ (obtained by numerically solving the traveling wave equation~\eqref{eq:TW_ode}) extended trivially in the $y$-direction. The PDE~\eqref{eq:grda_mk} was solved in a co-moving frame; finite differences were used for spatial discretization and Matlab's ode15s routine was used for time integration. The panels depict the $U$-profile at times $t=\{10000, 30000, 60000, 90000\}$ (from left to right). Yello indicates the bare soil state, while green indicates a vegetated state.  We observe that after long simulation times, the front interface exhibits bounded cusping behavior, reminiscent of Kuramoto-Sivashinsky dynamics associated to a sideband instability~\cite{doelman2009dynamics,hyman1986kuramoto}. While the same cusping behavior appeared for these parameters at smaller values of $\nu$, finger-like protrusions may appear for smaller values of $\nu$ in other parameter regimes; see Figure~\ref{fig:finger}.}
\label{fig:cusping}
\end{figure}

\begin{figure}
\hspace{.01\textwidth}
\begin{subfigure}{.22 \textwidth}
\centering
\includegraphics[width=1\linewidth]{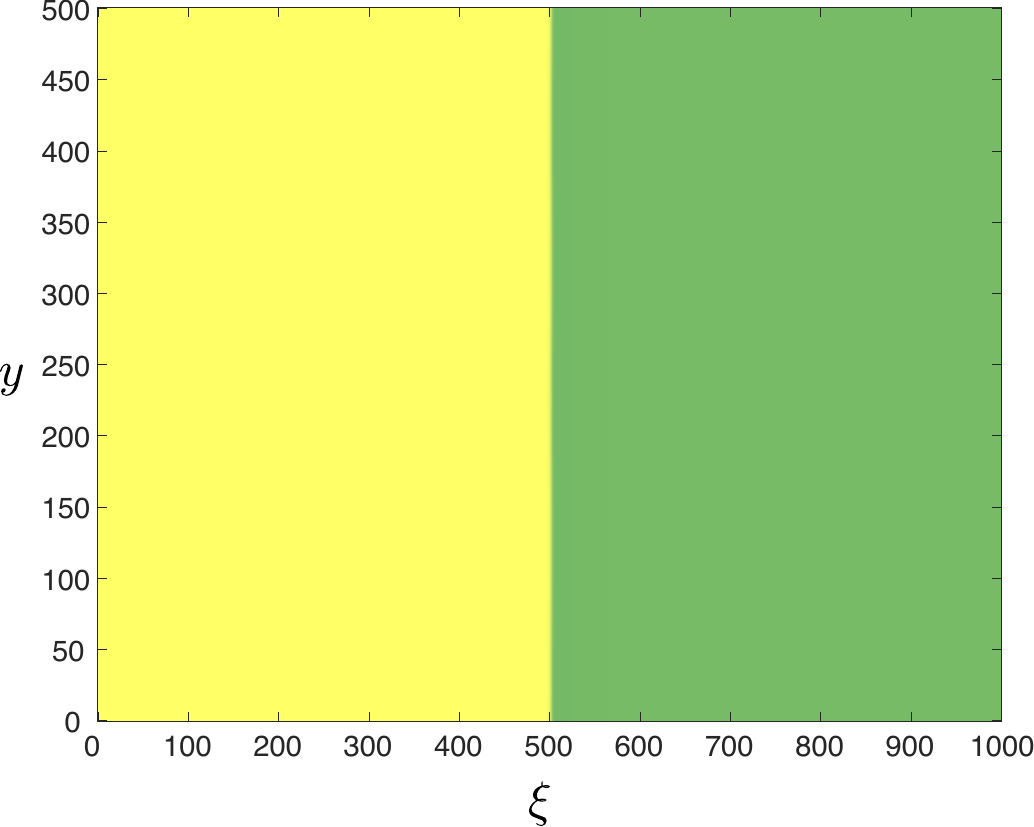}
\end{subfigure}
\hspace{.01\textwidth}
\begin{subfigure}{.22 \textwidth}
\centering
\includegraphics[width=1\linewidth]{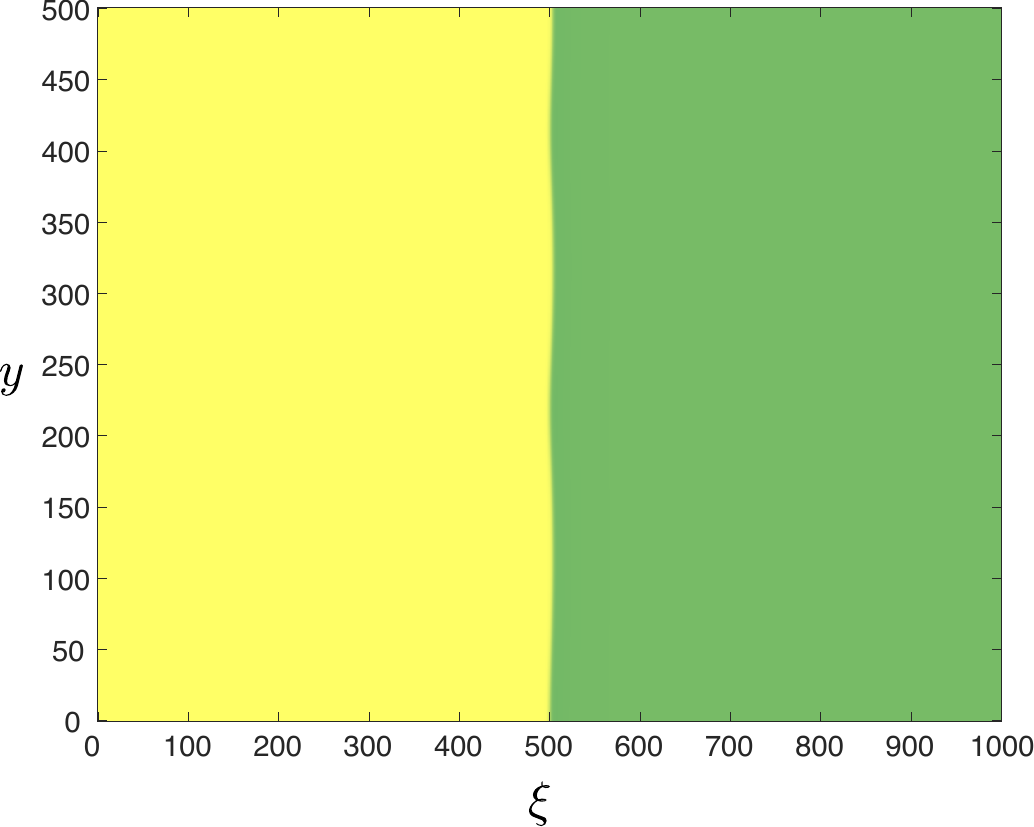}
\end{subfigure}
\hspace{.01\textwidth}
\begin{subfigure}{.22 \textwidth}
\centering
\includegraphics[width=1\linewidth]{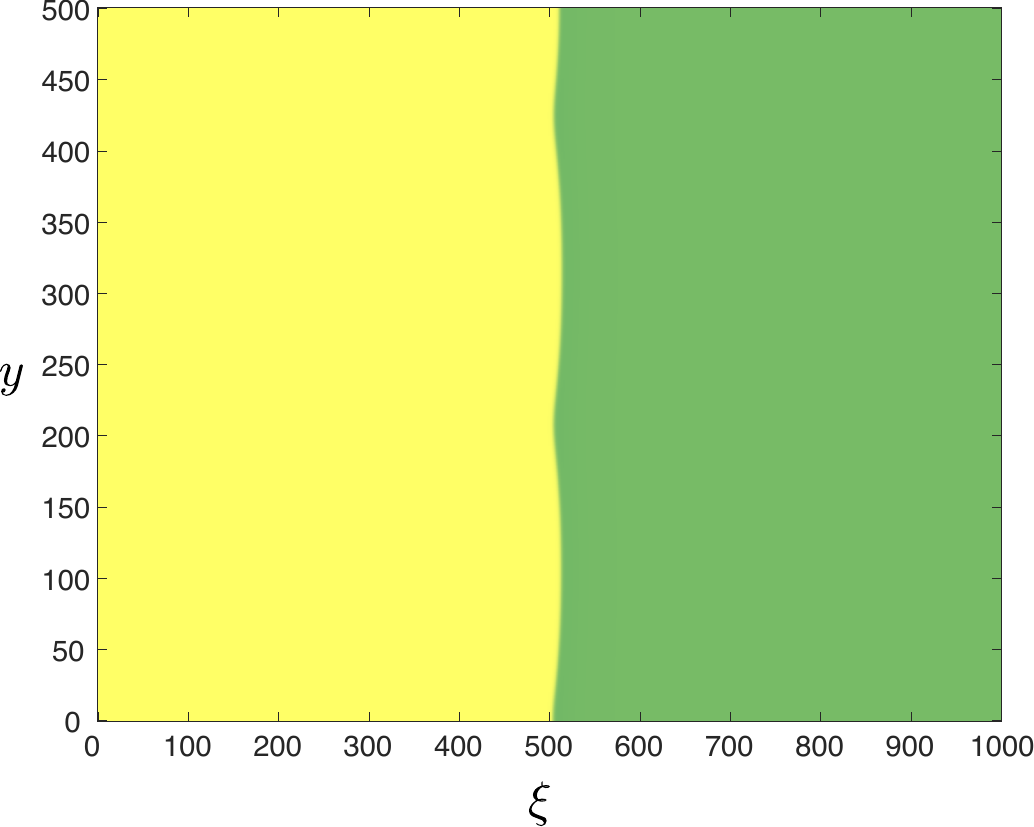}
\end{subfigure}
\hspace{.01\textwidth}
\begin{subfigure}{.22 \textwidth}
\centering
\includegraphics[width=1\linewidth]{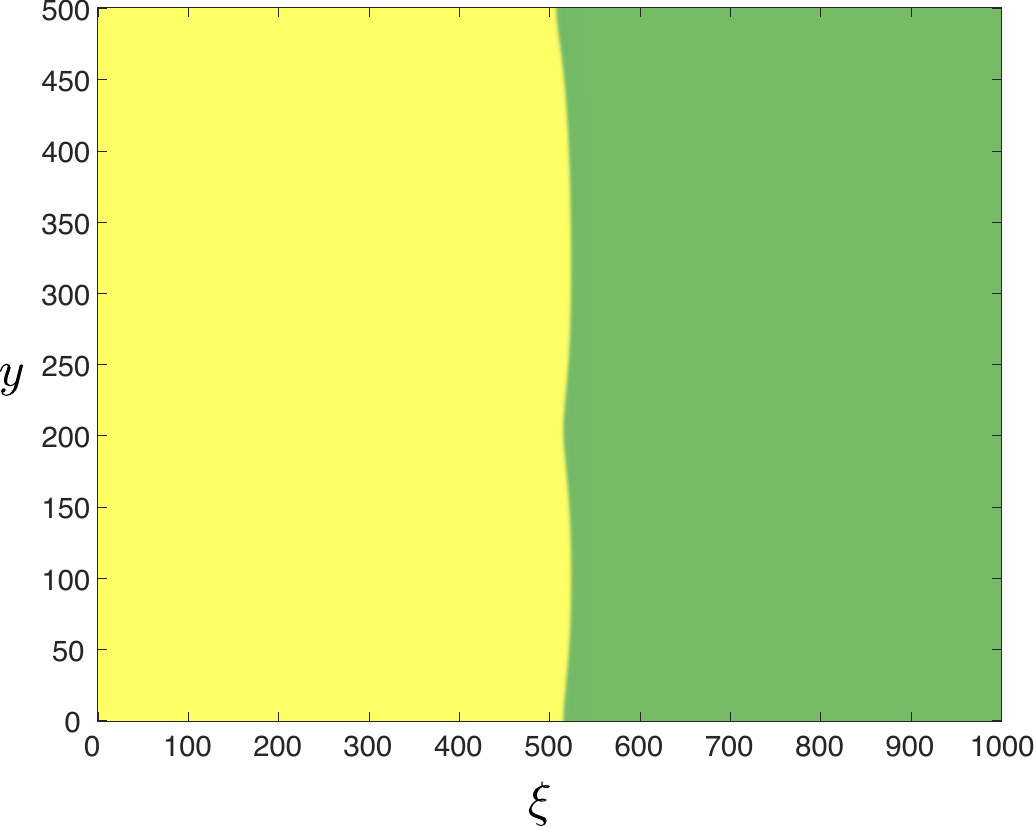}
\end{subfigure}\\

\hspace{.01\textwidth}
\begin{subfigure}{.22 \textwidth}
\centering
\includegraphics[width=1\linewidth]{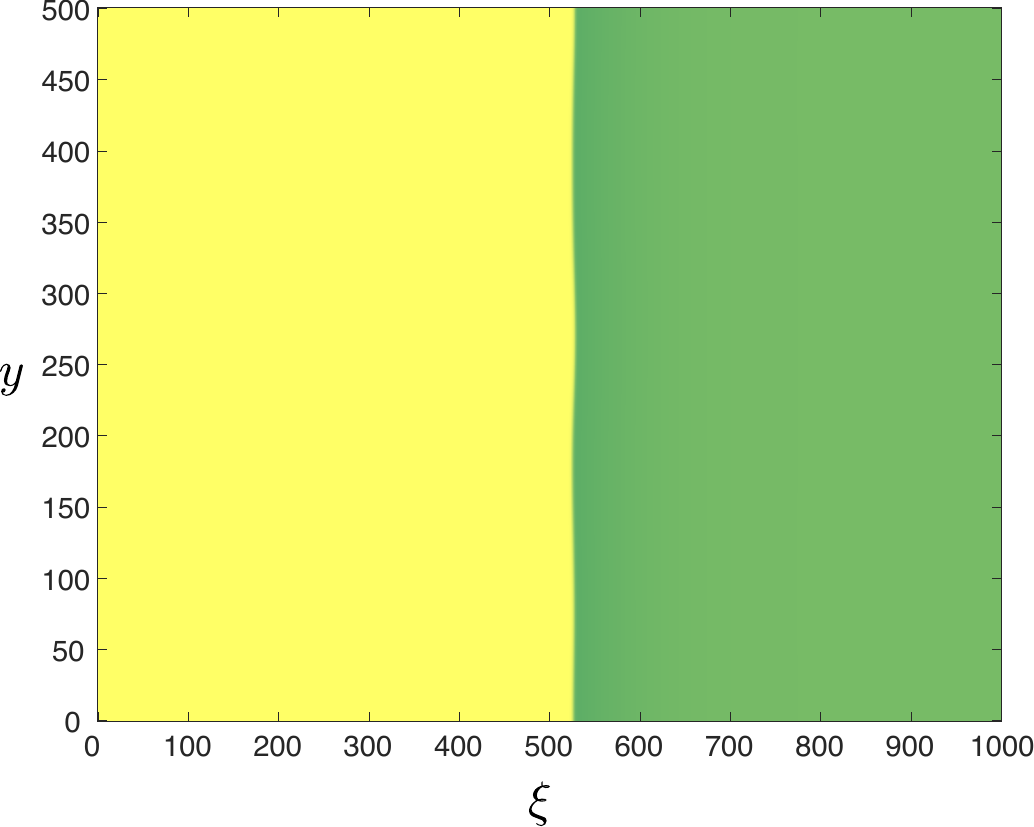}
\end{subfigure}
\hspace{.01\textwidth}
\begin{subfigure}{.22 \textwidth}
\centering
\includegraphics[width=1\linewidth]{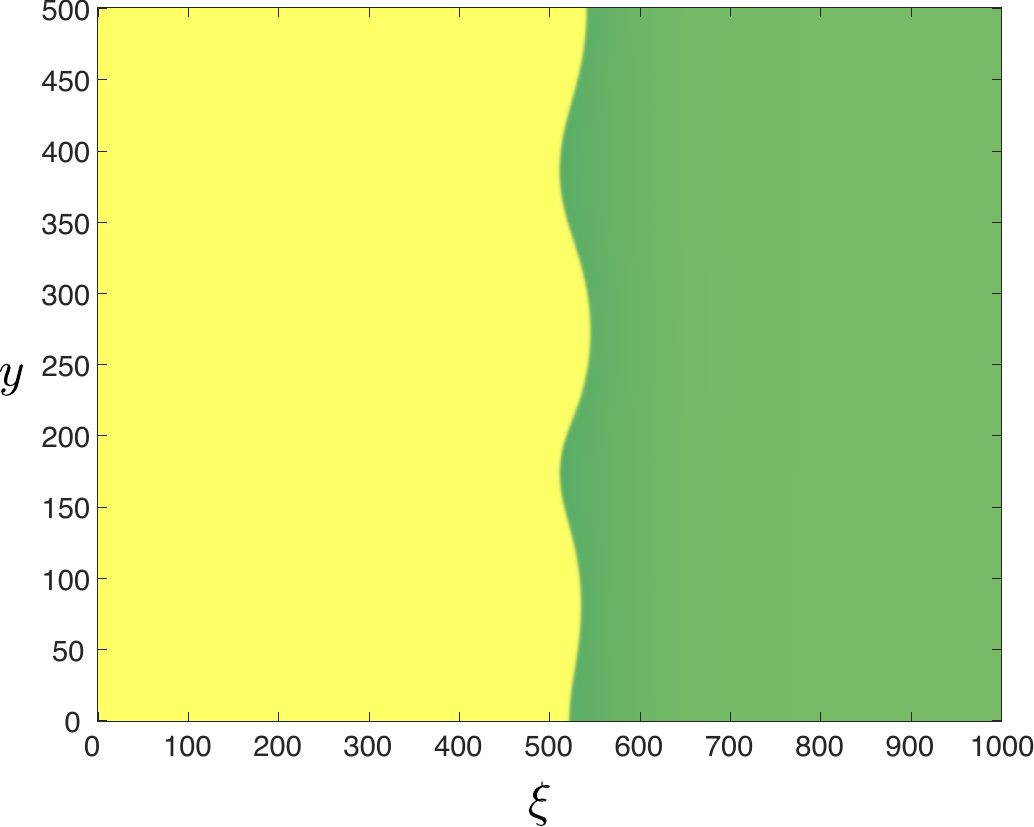}
\end{subfigure}
\hspace{.01\textwidth}
\begin{subfigure}{.22 \textwidth}
\centering
\includegraphics[width=1\linewidth]{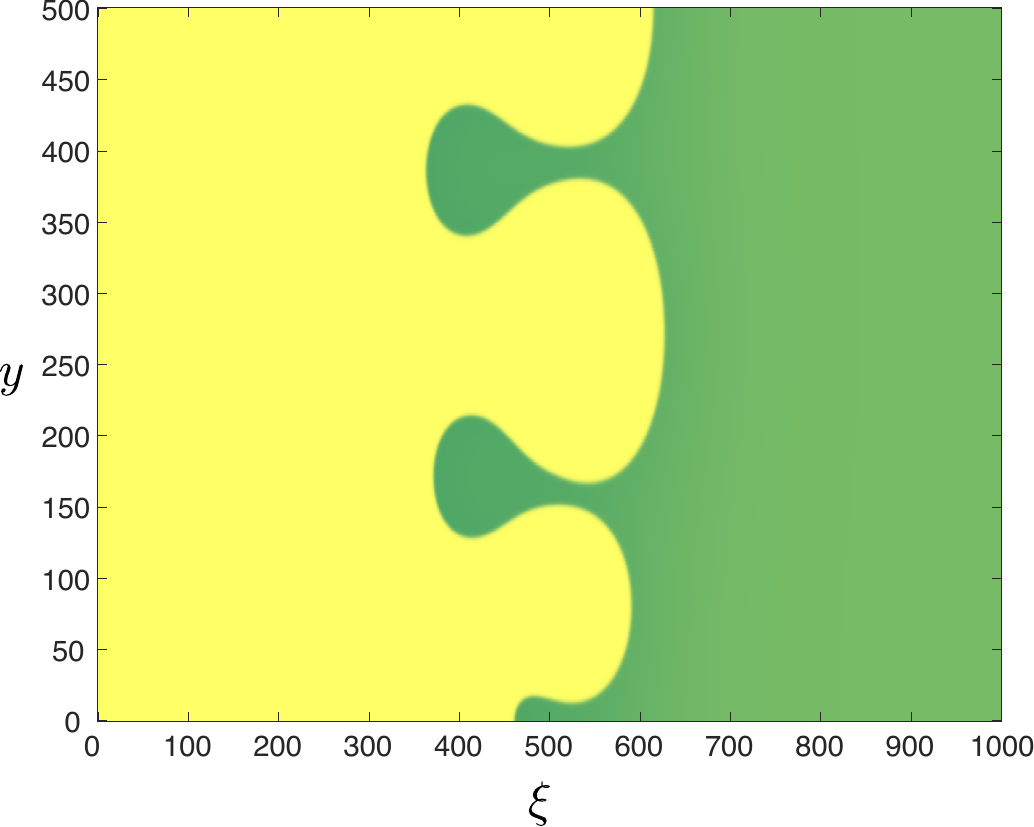}
\end{subfigure}
\hspace{.01\textwidth}
\begin{subfigure}{.22 \textwidth}
\centering
\includegraphics[width=1\linewidth]{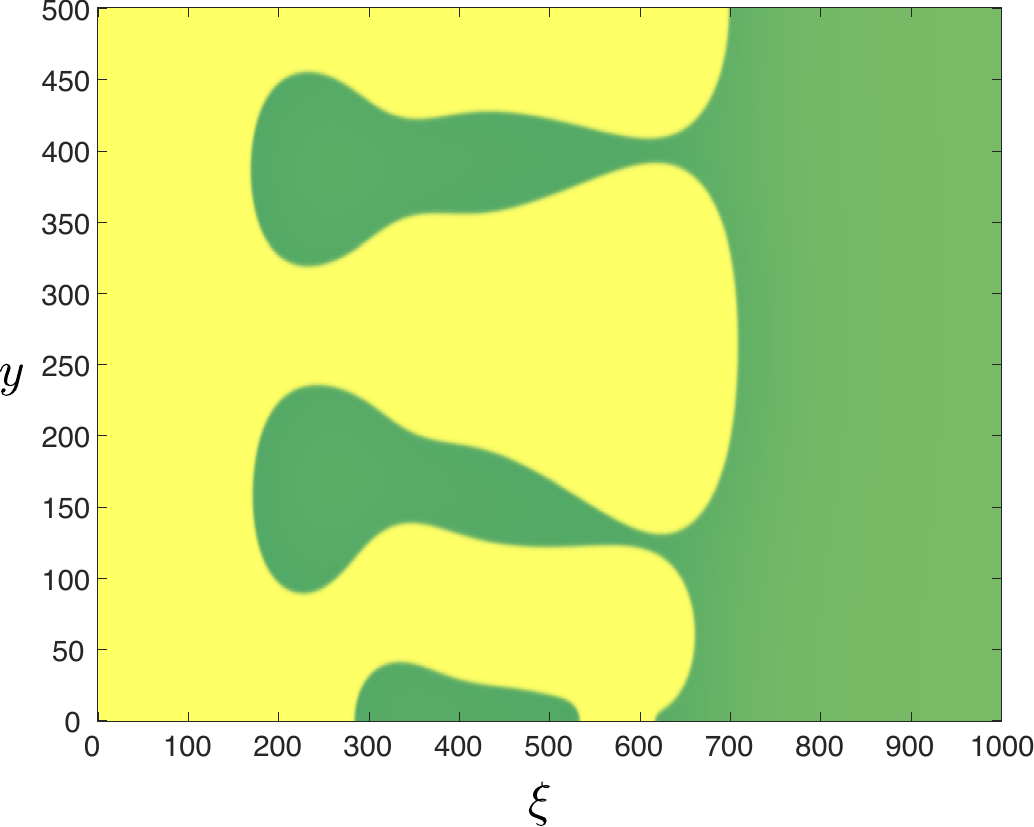}
\end{subfigure}
\caption{Shown are the results of direct numerical simulations in~\eqref{eq:grda_mk} for $\delta=0.01$ and $\boldsymbol{\mu}=(1.2,1.0,6.2)$, performed similarly as those in Figure~\ref{fig:cusping}. The top row depicts the $U$-profile of the solution in the case $\nu=500$ at times $t=\{5000, 10000, 15000, 20000\}$, while the bottom row depicts the solution in the case $\nu=50$ at times $t=\{2000, 3000, 4000, 5000\}$.  The critical stability boundary at which $\lambda_{\mathrm{c},2}$ changes sign was estimated at $\nu\approx 819$ for these parameter values. In the case $\nu=500$ (top row) closer to the stability boundary, the same cusping behavior is observed as in Figure~\ref{fig:cusping}, while in the case $\nu=50$ (bottom row) finger-like protrusions emerge on a comparatively shorter timescale.}
\label{fig:finger}
\end{figure}

These initial findings suggest numerous avenues for investigating the interplay of advection and diffusion on stability problems for traveling waves beyond one spatial dimension. In particular, for systems such as~\eqref{eq:modifiedKlausmeier}, where the quantities $F_*, G_*$ are such that $\lambda_{\mathrm{c},2}>0$ when $\nu=0$, we expect an exchange of stability upon increasing the advection $\nu$, occurring at a critical scaling $\nu\sim \delta^{-4/3}$. However, the coefficient $\lambda_{\mathrm{c},2}$ gives only a spectral stability criterion, so a natural question concerns the manifestation of the instability in the full nonlinear dynamics of~\eqref{eq:grda} which manifests as $\nu$ decreases through this critical value. Such interfacial instabilities (in the absence of advection) have been examined in other ecosystems models, see e.g.~\cite{CDLOR,fernandez2019front}, where the emergence of finger-like protrusions and cusps have been observed. Figures~\ref{fig:cusping}-\ref{fig:finger} show the result of direct numerical simulations in the model~\eqref{eq:grda_mk}. We observe the appearance of bounded cusp-like instabilities appearing in the planar interface for values of $\nu$ near (but below) the critical scaling $\nu\sim \delta^{-4/3}$. While this cusping behavior also appears for smaller values of $\nu$ in some parameter regimes, in other regimes (see Figure~\ref{fig:finger}) we observe the development of finger-like instabilities in the interface, suggesting that this could perhaps be a more severe manifestation of this interfacial instability. An unfolding of the stability boundary at the critical scaling $\nu\sim \delta^{-4/3}$, and an analytical treatment of the emergent nonlinear behavior is the subject of ongoing work.

While the results presented here are of a formal asymptotic nature, we emphasize that they could be obtained rigorously using geometric singular perturbation techniques. In particular, in the strong advection regime $\nu\geq \mathcal{O}(\delta^{-2})$ we expect it is possible to show rigorously that the fronts described here (as well as multi-front stripe and periodic patterned solutions) are stable in two spatial dimensions by analyzing~\eqref{eq:mKstability} for all values of $\ell\in\mathbb{R}$ using a nearly identical approach as in~\cite[\S5]{BCD}, via exponential dichotomies and the Lin--Sandstede method~\cite{cdrs, lin1990using, sandstede1998stability}. A rigorous treatment of the stability of such structures is the subject of future work.






\paragraph{Acknowledgements:} The author gratefully acknowledges support from the National Science Foundation through the grant DMS-2105816.

\appendix

\section{Stability of steady states}\label{app:steadystates}
We briefly derive the conditions~\eqref{eq:steadystate_conditions} of Assumption~\ref{assump:steadystates}. Linearizing about the steady states $(U^\pm,V^\pm)$ of~\eqref{eq:grda} and setting $(U,V)(x,y,t)=(U^\pm,V^\pm)+e^{\lambda t+ikx+i\ell y}(\bar{U}, \bar{V})$, we obtain the eigenvalue problem
\begin{align}\label{eq:steadystate_eigenvalueproblem}
\lambda \begin{pmatrix} \bar{U}\\ \bar{V}\end{pmatrix} = \begin{pmatrix}-(k^2+\ell^2) +F_u^\pm & F_v^\pm\\ G_u^\pm& -\frac{1}{\delta^2}(k^2+\ell^2)+i\nu k+ G_v^\pm\end{pmatrix} \begin{pmatrix} \bar{U}\\ \bar{V}\end{pmatrix}.
\end{align}
The steady states $(U^\pm,V^\pm)$ are stable provided $\Re \lambda(k,\ell)<0$ for any eigenvalue $\lambda(k,\ell)$ of~\eqref{eq:steadystate_eigenvalueproblem}, or equivalently, any $\lambda(k,\ell)$ satisfying
\begin{align}
    \lambda^2+(p_1+iq_1)\lambda  +(p_2+iq_2)=0
\end{align}
where
\begin{align*}
p_1&= \left(1+\frac{1}{\delta^2}\right)(k^2+\ell^2) - F_u^\pm-G_v^\pm\\
q_1&=- \nu k\\
p_2 &= F_u^\pm G_v^\pm-F_v^\pm G_u^\pm-\left(\frac{F_u^\pm}{\delta^2}+G_v^\pm\right) (k^2+\ell^2) +\frac{1}{\delta^2}(k^2+\ell^2)^2\\
q_2&= -\nu k(k^2+\ell^2)+ \nu k F_u^\pm.
\end{align*}
The corresponding roots of this quadratic have strictly negative real part if and only if (see e.g.~\cite{frank1946zeros})
\begin{align}\label{eq:app_condi}
    p_1&>0\\
    p_1^2p_2+p_1q_1q_2-q_2^2&>0 \label{eq:app_condii}
\end{align}
are satisfied for all $k, \ell \in \mathbb{R}$. The condition~\eqref{eq:app_condi} implies that
\begin{align}
    F_u^\pm+G_v^\pm<0
\end{align}
while~\eqref{eq:app_condii} implies (consider, e.g. $k=\ell=0$)
\begin{align}
    F_u^\pm G_v^\pm-F_v^\pm G_u^\pm>0.
\end{align}
Considering~\eqref{eq:app_condii} for $k=0$ and $\ell=\mathcal{O}(\delta)$ implies that $F_u^\pm<0$, while setting $\ell=0, k\ll1$ and $\nu$ sufficiently large implies that $G_v^\pm<0$ as well.

\bibliographystyle{abbrv}
\bibliography{bibfile}

\begin{thebibliography}{10}

\bibitem{banerjee2023rethinking}
S.~Banerjee, M.~Baudena, P.~Carter, R.~Bastiaansen, A.~Doelman, and
  M.~Rietkerk.
\newblock Rethinking tipping points in spatial ecosystems.
\newblock {\em arXiv preprint arXiv:2306.13571}, 2023.

\bibitem{barbier2014case}
N.~Barbier, P.~Couteron, and V.~Deblauwe.
\newblock Case study of self-organized vegetation patterning in dryland regions
  of central {A}frica.
\newblock In {\em Patterns of land degradation in drylands}, pages 347--356.
  Springer, 2014.

\bibitem{BCD}
R.~Bastiaansen, P.~Carter, and A.~Doelman.
\newblock Stable planar vegetation stripe patterns on sloped terrain in dryland
  ecosystems.
\newblock {\em Nonlinearity}, 32(8):2759, 2019.

\bibitem{bastiaansen2022fragmented}
R.~Bastiaansen, H.~A. Dijkstra, and A.~S. von~der Heydt.
\newblock Fragmented tipping in a spatially heterogeneous world.
\newblock {\em Environmental Research Letters}, 17(4):045006, 2022.

\bibitem{bennett2019large}
J.~J. Bennett and J.~A. Sherratt.
\newblock Large scale patterns in mussel beds: stripes or spots?
\newblock {\em Journal of mathematical biology}, 78(3):815--835, 2019.

\bibitem{borgogno2009mathematical}
F.~Borgogno, P.~D'Odorico, F.~Laio, and L.~Ridolfi.
\newblock Mathematical models of vegetation pattern formation in ecohydrology.
\newblock {\em Reviews of geophysics}, 47(1), 2009.

\bibitem{borthagaray2010vegetation}
A.~I. Borthagaray, M.~A. Fuentes, and P.~A. Marquet.
\newblock Vegetation pattern formation in a fog-dependent ecosystem.
\newblock {\em Journal of Theoretical Biology}, 265(1):18--26, 2010.

\bibitem{BCDL}
E.~Byrnes, P.~Carter, A.~Doelman, and L.~Liu.
\newblock Large amplitude radially symmetric spots and gaps in a dryland
  ecosystem model.
\newblock {\em Journal of Nonlinear Science}, 33(6):107, 2023.

\bibitem{cdrs}
P.~Carter, B.~de~Rijk, and B.~Sandstede.
\newblock Stability of traveling pulses with oscillatory tails in the
  {F}itz{H}ugh--{N}agumo system.
\newblock {\em Journal of Nonlinear Science}, 26(5):1369--1444, 2016.

\bibitem{CDLOR}
P.~Carter, A.~Doelman, K.~Lilly, E.~Obermayer, and S.~Rao.
\newblock Criteria for the (in)stability of planar interfaces in singularly
  perturbed 2-component reaction--diffusion equations.
\newblock {\em Physica D: Nonlinear Phenomena}, 444:133596, 2023.

\bibitem{chen2008evolution}
X.~Chen, R.~Hambrock, and Y.~Lou.
\newblock Evolution of conditional dispersal: a reaction--diffusion--advection
  model.
\newblock {\em Journal of mathematical biology}, 57(3):361--386, 2008.

\bibitem{deblauwe2012determinants}
V.~Deblauwe, P.~Couteron, J.~Bogaert, and N.~Barbier.
\newblock Determinants and dynamics of banded vegetation pattern migration in
  arid climates.
\newblock {\em Ecological monographs}, 82(1):3--21, 2012.

\bibitem{deblauwe2011environmental}
V.~Deblauwe, P.~Couteron, O.~Lejeune, J.~Bogaert, and N.~Barbier.
\newblock Environmental modulation of self-organized periodic vegetation
  patterns in {S}udan.
\newblock {\em Ecography}, 34(6):990--1001, 2011.

\bibitem{doelman2022slow}
A.~Doelman.
\newblock Slow localized patterns in singularly perturbed two-component
  reaction--diffusion equations.
\newblock {\em Nonlinearity}, 35(7):3487, 2022.

\bibitem{doelman2009dynamics}
A.~Doelman, B.~Sandstede, A.~Scheel, and G.~Schneider.
\newblock {\em The dynamics of modulated wave trains}.
\newblock American Mathematical Soc., 2009.

\bibitem{eigentler2021species}
L.~Eigentler.
\newblock Species coexistence in resource-limited patterned ecosystems is
  facilitated by the interplay of spatial self-organisation and intraspecific
  competition.
\newblock {\em Oikos}, 130(4):609--623, 2021.

\bibitem{fenichel1979geometric}
N.~Fenichel.
\newblock Geometric singular perturbation theory for ordinary differential
  equations.
\newblock {\em Journal of differential equations}, 31(1):53--98, 1979.

\bibitem{fernandez2019front}
C.~Fernandez-Oto, O.~Tzuk, and E.~Meron.
\newblock Front instabilities can reverse desertification.
\newblock {\em Physical review letters}, 122(4):048101, 2019.

\bibitem{frank1946zeros}
E.~Frank.
\newblock On the zeros of polynomials with complex coefficients.
\newblock {\em Bulletin of the American Mathematical Society}, 52(2):144--157,
  1946.

\bibitem{gandhi2018influence}
P.~Gandhi, L.~Werner, S.~Iams, K.~Gowda, and M.~Silber.
\newblock A topographic mechanism for arcing of dryland vegetation bands.
\newblock {\em Journal of The Royal Society Interface}, 15(147), 2018.

\bibitem{ge2015sis}
J.~Ge, K.~I. Kim, Z.~Lin, and H.~Zhu.
\newblock A sis reaction--diffusion--advection model in a low-risk and
  high-risk domain.
\newblock {\em Journal of Differential Equations}, 259(10):5486--5509, 2015.

\bibitem{gilad2004ecosystem}
E.~Gilad, J.~von Hardenberg, A.~Provenzale, M.~Shachak, and E.~Meron.
\newblock Ecosystem engineers: from pattern formation to habitat creation.
\newblock {\em Physical Review Letters}, 93(9):098105, 2004.

\bibitem{hellden1988desertification}
U.~Helld{\'e}n.
\newblock Desertification monitoring: {I}s the desert encroaching?
\newblock {\em Desertification Control Bulletin}, 17:8--12, 1988.

\bibitem{hyman1986kuramoto}
J.~M. Hyman and B.~Nicolaenko.
\newblock The {K}uramoto-{S}ivashinsky equation: a bridge between {PDE}'s and
  dynamical systems.
\newblock {\em Physica D: Nonlinear Phenomena}, 18(1-3):113--126, 1986.

\bibitem{klausmeier1999regular}
C.~A. Klausmeier.
\newblock Regular and irregular patterns in semiarid vegetation.
\newblock {\em Science}, 284(5421):1826--1828, 1999.

\bibitem{kolokolnikov2018stabilizing}
T.~Kolokolnikov, M.~Ward, J.~Tzou, and J.~Wei.
\newblock Stabilizing a homoclinic stripe.
\newblock {\em Philosophical Transactions of the Royal Society A: Mathematical,
  Physical and Engineering Sciences}, 376(2135):20180110, 2018.

\bibitem{lin1990using}
X.-B. Lin.
\newblock Using melnikov's method to solve silnikov's problems.
\newblock {\em Proceedings of the Royal Society of Edinburgh Section A:
  Mathematics}, 116(3-4):295--325, 1990.

\bibitem{lipcius2021facilitation}
R.~N. Lipcius, D.~G. Matthews, L.~Shaw, J.~Shi, and S.~Zaytseva.
\newblock Facilitation between ecosystem engineers, salt marsh grass and
  mussels, produces pattern formation on salt marsh shorelines.
\newblock {\em bioRxiv}, 2021.

\bibitem{ludwig2005vegetation}
J.~Ludwig, B.~Wilcox, D.~Breshears, D.~Tongway, and A.~Imeson.
\newblock Vegetation patches and runoff-erosion as interacting ecohydrological
  processes in semiarid landscapes.
\newblock {\em Ecology}, 86(2):288--297, 2005.

\bibitem{malchow1996nonlinear}
H.~Malchow.
\newblock Nonlinear plankton dynamics and pattern formation in an
  ecohydrodynamic model system.
\newblock {\em Journal of Marine Systems}, 7(2-4):193--202, 1996.

\bibitem{meron2012pattern}
E.~Meron.
\newblock Pattern-formation approach to modelling spatially extended
  ecosystems.
\newblock {\em Ecological Modelling}, 234:70--82, 2012.

\bibitem{rietkerk2008regular}
M.~Rietkerk and J.~Van~de Koppel.
\newblock Regular pattern formation in real ecosystems.
\newblock {\em Trends in ecology \& evolution}, 23(3):169--175, 2008.

\bibitem{rovinsky1992chemical}
A.~B. Rovinsky and M.~Menzinger.
\newblock Chemical instability induced by a differential flow.
\newblock {\em Physical Review Letters}, 69(8):1193, 1992.

\bibitem{sandstede1998stability}
B.~Sandstede.
\newblock Stability of multiple-pulse solutions.
\newblock {\em Transactions of the American Mathematical Society},
  350(2):429--472, 1998.

\bibitem{sewaltspatially}
L.~Sewalt and A.~Doelman.
\newblock Spatially periodic multipulse patterns in a generalized
  {K}lausmeier--{G}ray--{S}cott model.
\newblock {\em SIAM Journal on Applied Dynamical Systems}, 16(2):1113--1163,
  2017.

\bibitem{siero2015striped}
E.~Siero, A.~Doelman, M.~Eppinga, J.~D. Rademacher, M.~Rietkerk, and K.~Siteur.
\newblock Striped pattern selection by advective reaction-diffusion systems:
  Resilience of banded vegetation on slopes.
\newblock {\em Chaos: An Interdisciplinary Journal of Nonlinear Science},
  25(3):036411, 2015.

\bibitem{siteur2014beyond}
K.~Siteur, E.~Siero, M.~B. Eppinga, J.~D. Rademacher, A.~Doelman, and
  M.~Rietkerk.
\newblock Beyond {T}uring: The response of patterned ecosystems to
  environmental change.
\newblock {\em Ecological Complexity}, 20:81--96, 2014.

\bibitem{taniguchi2003instability}
M.~Taniguchi.
\newblock Instability of planar traveling waves in bistable reaction-diffusion
  systems.
\newblock {\em Discrete \& Continuous Dynamical Systems-B}, 3(1):21, 2003.

\bibitem{taniguchi1994instability}
M.~Taniguchi and Y.~Nishiura.
\newblock Instability of planar interfaces in reaction-diffusion systems.
\newblock {\em SIAM Journal on Mathematical Analysis}, 25(1):99--134, 1994.

\bibitem{un2015transforming}
G.~A. UN.
\newblock Transforming our world: The 2030 agenda for sustainable development.
\newblock Technical report, A/RES/70/1, 21 October, 2015.

\bibitem{valentin1999soil}
C.~Valentin, J.~d'Herb{\`e}s, and J.~Poesen.
\newblock Soil and water components of banded vegetation patterns.
\newblock {\em Catena}, 37(1):1--24, 1999.

\bibitem{zelnik2017desertification}
Y.~R. Zelnik, H.~Uecker, U.~Feudel, and E.~Meron.
\newblock Desertification by front propagation?
\newblock {\em Journal of Theoretical Biology}, 418:27--35, 2017.

\end{thebibliography}

\end{document}